\documentclass[a4paper]{amsart}
\usepackage[breaklinks,pdfstartview=FitH,colorlinks,linktocpage,linkcolor=blue,citecolor=blue,urlcolor=blue]{hyperref}

\usepackage{amssymb}
\usepackage{amsthm}

\usepackage{upref}
\usepackage{geometry}

\usepackage{pgf,tikz}\usetikzlibrary{%
   arrows%
}

\newtheorem{theorem}{Theorem}
\newtheorem*{theorem*}{Theorem}
\newtheorem{claim}[theorem]{Claim}

\theoremstyle{definition}
\newtheorem*{definition*}{Definition}
\newtheorem{example}[theorem]{Example}
\newtheorem{problem}[theorem]{Problem}

\begin{document}

\definecolor{ffzzzz}{rgb}{1,0.6,0.6}
\definecolor{wwttqq}{rgb}{0.4,0.2,0}
\definecolor{tttttt}{rgb}{0.2,0.2,0.2}
\definecolor{wwwwww}{rgb}{0.4,0.4,0.4}
\definecolor{cccccc}{rgb}{0.8,0.8,0.8}
\definecolor{zzqqcc}{rgb}{0.6,0,0.8}
\definecolor{zzzzzz}{rgb}{0.6,0.6,0.6}

\large

\title{Packing a cake into a box}


\author{M. Skopenkov}
\address{Institute for information transmission problems of the Russian Academy of Sciences\\
Bolshoy Karetny per. 19, bld. 1, Moscow, 127994, Russian Federation\\
and King Abdullah University of Science and Technology\\
4700 Thuwal, 23955-6900, Kingdom of Saudi Arabia}
\email{skopenkov@rambler.ru}

\date{}

\begin{abstract}
Given a cake in form of a triangle and a box that fits the mirror image
of the cake, how to cut the cake into a minimal number of pieces so that it
can be put into the box? The cake has an icing, so that we are not
allowed to put it into the box upside down. V.G. Boltyansky asked
this question in 1977 and showed that three pieces always suffice.
In this paper we provide examples of cakes that cannot be cut into two
pieces to put into the box. This shows that three is the answer to
V.G. Boltyansky's question. Also we give examples of cakes
which can be cut into two pieces.
\end{abstract}






\maketitle

\section{Examples}

\begin{figure}[htb]
\includegraphics[width=8cm]{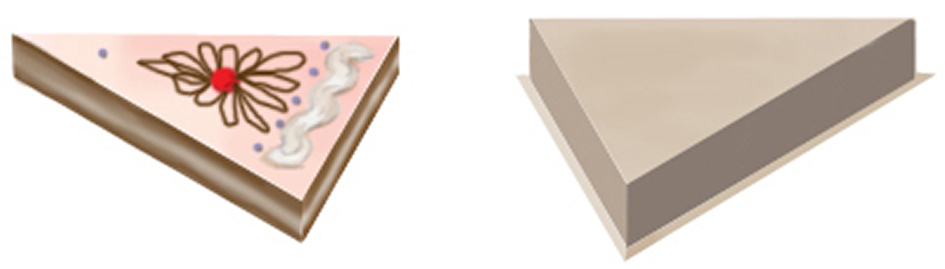}
\caption{A cake and a box}
\label{cake}
\end{figure}

Given a cake in form of a triangle and a box that fits the mirror image
of the cake, how to cut the cake into a minimal number of pieces so that it
can be put into the box? The cake has an icing, so that we are not
allowed to put it into the box upside down, see Figure~\ref{cake}.

\begin{figure}[htb]
{\fontsize{9pt}{9pt}\selectfont\begin{tikzpicture}[line cap=round,line join=round,>=triangle 45,x=1.0cm,y=1.0cm]
\clip(7.82,2.88) rectangle (15.82,4.86);
\fill[color=ffzzzz,fill=ffzzzz,fill opacity=0.2] (8.08,4.72) -- (10.48,3.34) -- (11.3,4.72) -- cycle;
\draw [shift={(8.08,4.72)},color=wwttqq,fill=wwttqq,fill opacity=0.2] (0,0) -- (-29.9:0.8) arc (-29.9:0:0.8) -- cycle;
\draw [shift={(10.48,3.34)},color=wwttqq,fill=wwttqq,fill opacity=0.2] (0,0) -- (59.28:0.6) arc (59.28:150.1:0.6) -- cycle;
\draw [shift={(11.3,4.72)},color=wwttqq,fill=wwttqq,fill opacity=0.2] (0,0) -- (180:0.7) arc (180:239.28:0.7) -- cycle;
\fill[color=ffzzzz,fill=ffzzzz,fill opacity=0.2] (15.65,4.72) -- (13.25,3.34) -- (12.43,4.72) -- cycle;
\draw [shift={(15.65,4.72)},color=wwttqq,fill=wwttqq,fill opacity=0.2] (0,0) -- (180:0.8) arc (180:209.9:0.8) -- cycle;
\draw [shift={(12.43,4.72)},color=wwttqq,fill=wwttqq,fill opacity=0.2] (0,0) -- (-59.28:0.7) arc (-59.28:0:0.7) -- cycle;
\draw [shift={(13.25,3.34)},color=wwttqq,fill=wwttqq,fill opacity=0.2] (0,0) -- (29.9:0.6) arc (29.9:120.72:0.6) -- cycle;
\draw [color=wwttqq] (8.08,4.72)-- (10.48,3.34);
\draw [color=wwttqq] (10.48,3.34)-- (11.3,4.72);
\draw [color=wwttqq] (11.3,4.72)-- (8.08,4.72);
\draw [dash pattern=on 1pt off 3pt on 5pt off 4pt,color=wwttqq] (11.86,2.88) -- (11.86,4.86);
\draw [color=wwttqq] (15.65,4.72)-- (13.25,3.34);
\draw [color=wwttqq] (13.25,3.34)-- (12.43,4.72);
\draw [color=wwttqq] (12.43,4.72)-- (15.65,4.72);
\draw [color=wwttqq](23.54,7.6) node[anchor=north west] {$\alpha$};
\draw (0,0.24) node[anchor=north west] {$\alpha$};
\fill [color=wwttqq] (8.08,4.72) circle (1.0pt);
\draw[color=wwttqq] (8.08,4.46) node {$A$};
\fill [color=wwttqq] (10.48,3.34) circle (1.0pt);
\draw[color=wwttqq] (10.64,3.12) node {$B$};
\fill [color=wwttqq] (11.3,4.72) circle (1.0pt);
\draw[color=wwttqq] (11.4,4.44) node {$C$};
\draw[color=wwttqq] (8.72,4.56) node {$\alpha$};
\draw[color=wwttqq] (10.4,3.66) node {$\beta$};
\draw[color=wwttqq] (10.94,4.48) node {$\gamma$};
\fill [color=wwttqq] (15.65,4.72) circle (1.0pt);
\draw[color=wwttqq] (15.52,4.46) node {$A'$};
\fill [color=wwttqq] (13.25,3.34) circle (1.0pt);
\draw[color=wwttqq] (13.12,3.14) node {$B'$};
\fill [color=wwttqq] (12.43,4.72) circle (1.0pt);
\draw[color=wwttqq] (12.32,4.44) node {$C'$};
\draw[color=wwttqq] (15.06,4.54) node {$\alpha$};
\draw[color=wwttqq] (12.86,4.48) node {$\gamma$};
\draw[color=wwttqq] (13.34,3.64) node {$\beta$};
\end{tikzpicture}}
\caption{A cake and a box again}
\label{fig1}
\end{figure}
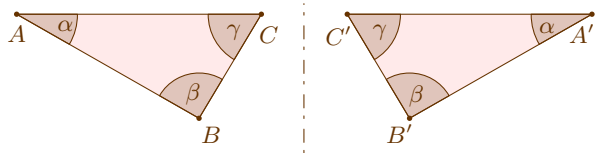

V.G. Boltyansky asked this question in his book on Hilbert's third problem \cite[\S9]{Bol77}. 
Let us give an accurate statement, see Figure~\ref{fig1}. 
A \emph{cake}
\begin{tikzpicture}[line cap=round,line join=round,>=triangle 45,x=0.11452372190810167cm,y=0.11452372190810167cm]
\clip(-1.75,3.88) rectangle (1.74,5.89);
\draw (-1.72,5.78)-- (0.72,4.38);
\draw (0.72,4.38)-- (1.56,5.78);
\draw (-1.72,5.78)-- (1.56,5.78);
\draw (-1.72,5.78)-- (-1.7,5.46);
\draw (-1.7,5.46)-- (0.72,4.08);
\draw (0.72,4.38)-- (0.72,4.08);
\draw (0.72,4.08)-- (1.56,5.46);
\draw (1.56,5.78)-- (1.56,5.46);
\draw [rotate around={-17.53:(0.13,5.38)}] (0.13,5.38) ellipse (0.03cm and 0.03cm);
\draw (0.42,4.8)-- (0.74,4.8);
\draw (0.74,4.8)-- (0.58,5.04);
\draw (0.58,5.04)-- (0.96,5.1);
\draw (0.96,5.1)-- (0.78,5.32);
\draw (0.78,5.32)-- (1.14,5.42);
\draw (1.14,5.42)-- (0.92,5.56);
\draw (0.92,5.56)-- (1.26,5.64);
\end{tikzpicture}
is a nonisosceles triangle in the plane with vertices $A$, $B$, $C$ and angles $\alpha$, $\beta$, $\gamma$. A \emph{box}
\begin{tikzpicture}[line cap=round,line join=round,>=triangle 45,x=0.11711167869707381cm,y=0.11711167869707381cm]
\clip(2.62,3.92) rectangle (6.04,6.09);
\draw (2.66,5.92)-- (5.92,5.92);
\draw (2.66,5.92)-- (3.48,4.52);
\draw (5.92,5.92)-- (3.48,4.52);
\draw (2.66,5.92)-- (2.66,5.44);
\draw (3.48,4.52)-- (3.5,4.02);
\draw (5.92,5.92)-- (5.94,5.46);
\draw (2.66,5.44)-- (3.5,4.02);
\draw (3.5,4.02)-- (5.94,5.46);
\draw (2.95,5.43)-- (5.09,5.44);
\end{tikzpicture}
is a triangle $A'B'C'$ obtained from the cake by a reflection with respect to a line. A cutting of the cake is \emph{nice} if the pieces are polygons which can be put into the box all together by means of appropriate rotations and translations.
The problem we study is to cut the cake nicely into a minimal number of pieces.

V.G. Boltyansky gives the following example:

\begin{example} \label{ex0}
\textup{\cite[\S9]{Bol77}}
Any cake can be nicely cut into $\mathbf{3}$ pieces.
\end{example}

\begin{figure}[htb]
\begin{tikzpicture}[line cap=round,line join=round,>=triangle 45,x=1.0cm,y=1.0cm]
\clip(7.86,3.26) rectangle (15.86,4.76);
\fill[color=wwttqq,fill=wwttqq,fill opacity=0.4] (8.08,4.72) -- (10.27,4.72) -- (10.27,4.13) -- (9.98,3.63) -- cycle;
\fill[color=ffzzzz,fill=ffzzzz,fill opacity=0.2] (11.3,4.72) -- (10.77,3.84) -- (10.27,4.13) -- (10.27,4.72) -- cycle;
\fill[color=zzzzzz,fill=zzzzzz,fill opacity=0.05] (10.48,3.34) -- (9.98,3.63) -- (10.27,4.13) -- (10.77,3.84) -- cycle;
\draw[dotted,color=zzqqcc] (10.27,4.51) -- (10.48,4.51) -- (10.48,4.72) -- (10.27,4.72) -- cycle;
\draw[dotted,color=zzqqcc] (10.59,3.94) -- (10.48,3.76) -- (10.67,3.65) -- (10.77,3.84) -- cycle;
\draw[dotted,color=zzqqcc] (10.09,3.81) -- (9.9,3.92) -- (9.8,3.73) -- (9.98,3.63) -- cycle;
\fill[color=wwttqq,fill=wwttqq,fill opacity=0.4] (15.65,4.72) -- (13.46,4.72) -- (13.46,4.13) -- (13.75,3.63) -- cycle;
\fill[color=ffzzzz,fill=ffzzzz,fill opacity=0.2] (12.43,4.72) -- (12.95,3.84) -- (13.46,4.13) -- (13.46,4.72) -- cycle;
\fill[color=zzzzzz,fill=zzzzzz,fill opacity=0.05] (13.25,3.34) -- (13.75,3.63) -- (13.46,4.13) -- (12.95,3.84) -- cycle;
\draw [color=wwttqq](23.54,7.6) node[anchor=north west] {$\alpha$};
\draw (0,0.24) node[anchor=north west] {$\alpha$};
\draw [color=wwttqq] (8.08,4.72)-- (10.27,4.72);
\draw [color=zzqqcc] (10.27,4.72)-- (10.27,4.13);
\draw [color=zzqqcc] (10.27,4.13)-- (9.98,3.63);
\draw [color=wwttqq] (9.98,3.63)-- (8.08,4.72);
\draw [color=ffzzzz] (11.3,4.72)-- (10.77,3.84);
\draw [color=zzqqcc] (10.77,3.84)-- (10.27,4.13);
\draw [color=zzqqcc] (10.27,4.13)-- (10.27,4.72);
\draw [color=ffzzzz] (10.27,4.72)-- (11.3,4.72);
\draw [color=zzzzzz] (10.48,3.34)-- (9.98,3.63);
\draw [color=zzqqcc] (9.98,3.63)-- (10.27,4.13);
\draw [color=zzqqcc] (10.27,4.13)-- (10.77,3.84);
\draw [color=zzzzzz] (10.77,3.84)-- (10.48,3.34);
\draw [dotted,color=zzqqcc] (10.27,4.13) circle (0.59cm);
\draw [->,color=ffzzzz] (10.53,4.43) -- (10.48,4.5);
\draw [->,color=wwttqq] (10.1,4.49) -- (10.16,4.54);
\draw [->,color=zzzzzz] (10.3,3.74) -- (10.06,3.76);
\draw [shift={(10.27,4.13)},color=wwttqq]  plot[domain=2.02:3.33,variable=\t]({1*0.39*cos(\t r)+0*0.39*sin(\t r)},{0*0.39*cos(\t r)+1*0.39*sin(\t r)});
\draw [shift={(10.27,4.13)},color=zzzzzz]  plot[domain=4.79:5.18,variable=\t]({1*0.39*cos(\t r)+0*0.39*sin(\t r)},{0*0.39*cos(\t r)+1*0.39*sin(\t r)});
\draw [shift={(10.27,4.13)},color=ffzzzz]  plot[domain=-0.22:0.86,variable=\t]({1*0.39*cos(\t r)+0*0.39*sin(\t r)},{0*0.39*cos(\t r)+1*0.39*sin(\t r)});
\draw [color=wwttqq] (15.65,4.72)-- (13.46,4.72);
\draw [color=zzqqcc] (13.46,4.72)-- (13.46,4.13);
\draw [color=zzqqcc] (13.46,4.13)-- (13.75,3.63);
\draw [color=wwttqq] (13.75,3.63)-- (15.65,4.72);
\draw [color=ffzzzz] (12.43,4.72)-- (12.95,3.84);
\draw [color=zzqqcc] (12.95,3.84)-- (13.46,4.13);
\draw [color=zzqqcc] (13.46,4.13)-- (13.46,4.72);
\draw [color=ffzzzz] (13.46,4.72)-- (12.43,4.72);
\draw [color=zzzzzz] (13.25,3.34)-- (13.75,3.63);
\draw [color=zzqqcc] (13.75,3.63)-- (13.46,4.13);
\draw [color=zzqqcc] (13.46,4.13)-- (12.95,3.84);
\draw [color=zzzzzz] (12.95,3.84)-- (13.25,3.34);
\fill [color=zzqqcc] (10.27,4.13) circle (1.0pt);
\fill [color=zzqqcc] (13.46,4.13) circle (1.0pt);
\fill [color=zzqqcc] (13.46,4.13) circle (1.0pt);
\fill [color=zzqqcc] (13.46,4.13) circle (1.0pt);
\end{tikzpicture}
\caption{A nice cutting into $3$ pieces}
\label{fig3parts}
\end{figure}
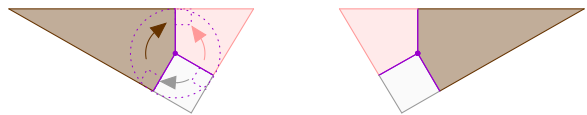

\begin{proof}
Indeed, cut the cake along the perpendiculars from the incenter to the sides, see Figure~\ref{fig3parts}.
Rotating the pieces around the incenter pack the cake into the box.
\end{proof}

This is all concerning the question we find in the book. However, V.V.~Prasolov indicated that some cakes
(see Example~\ref{ex1}(a)--(d)) can be nicely cut into $\mathbf{2}$ pieces.

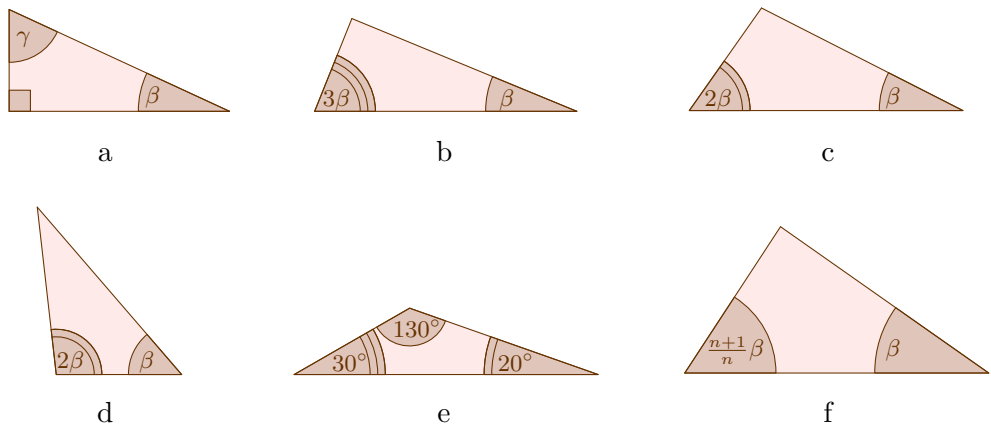
\begin{figure}[htb]
\begin{tabular}{ccc}
{\fontsize{9pt}{9pt}\selectfont\begin{tikzpicture}[line cap=round,line join=round,>=triangle 45,x=1.0cm,y=1.0cm]
\clip(-0.5,-0.08) rectangle (3.02,1.44);
\fill[color=ffzzzz,fill=ffzzzz,fill opacity=0.2] (0.05,0.05) -- (2.95,0.05) -- (0.05,1.4) -- cycle;
\draw[color=wwttqq,fill=wwttqq,fill opacity=0.2] (0.33,0.05) -- (0.33,0.33) -- (0.05,0.33) -- (0.05,0.05) -- cycle;
\draw [shift={(2.95,0.05)},color=wwttqq,fill=wwttqq,fill opacity=0.2] (0,0) -- (155.04:1.2) arc (155.04:180:1.2) -- cycle;
\draw [shift={(0.05,1.4)},color=wwttqq,fill=wwttqq,fill opacity=0.2] (0,0) -- (-90:0.7) arc (-90:-24.96:0.7) -- cycle;
\draw [color=wwttqq] (0.05,0.05)-- (2.95,0.05);
\draw [color=wwttqq] (2.95,0.05)-- (0.05,1.4);
\draw [color=wwttqq] (0.05,1.4)-- (0.05,0.05);
\draw[color=wwttqq] (1.94,0.24) node {$\beta$};
\draw[color=wwttqq] (0.24,1) node {$\gamma$};
\end{tikzpicture}}
&
{\fontsize{9pt}{9pt}\selectfont\begin{tikzpicture}[line cap=round,line join=round,>=triangle 45,x=1.0cm,y=1.0cm]
\clip(-0.3,-0.08) rectangle (3.7,1.42);
\fill[color=ffzzzz,fill=ffzzzz,fill opacity=0.2] (0.05,0.05) -- (3.5,0.05) -- (0.54,1.28) -- cycle;
\draw [shift={(0.05,0.05)},color=wwttqq,fill=wwttqq,fill opacity=0.2] (0,0) -- (0:0.8) arc (0:68.28:0.8) -- cycle;
\draw [shift={(3.5,0.05)},color=wwttqq,fill=wwttqq,fill opacity=0.2] (0,0) -- (157.44:1.2) arc (157.44:180:1.2) -- cycle;
\draw [color=wwttqq] (0.05,0.05)-- (3.5,0.05);
\draw [color=wwttqq] (3.5,0.05)-- (0.54,1.28);
\draw [color=wwttqq] (0.54,1.28)-- (0.05,0.05);
\draw [shift={(0.05,0.05)},color=wwttqq] (0:0.8) arc (0:68.28:0.8);
\draw [shift={(0.05,0.05)},color=wwttqq] (0:0.7) arc (0:68.28:0.7);
\draw [shift={(0.05,0.05)},color=wwttqq] (0:0.6) arc (0:68.28:0.6);
\draw[color=wwttqq] (0.34,0.2) node {$3\beta$};
\draw[color=wwttqq] (2.58,0.2) node {$\beta$};
\end{tikzpicture}}
&
{\fontsize{9pt}{9pt}\selectfont\begin{tikzpicture}[line cap=round,line join=round,>=triangle 45,x=1.0cm,y=1.0cm]
\clip(-0.22,-0.14) rectangle (3.76,1.5);
\fill[color=ffzzzz,fill=ffzzzz,fill opacity=0.2] (0,0) -- (3.6,0) -- (0.95,1.36) -- cycle;
\draw [shift={(0,0)},color=wwttqq,fill=wwttqq,fill opacity=0.2] (0,0) -- (0:0.8) arc (0:55.06:0.8) -- cycle;
\draw [shift={(3.6,0)},color=wwttqq,fill=wwttqq,fill opacity=0.2] (0,0) -- (152.83:1.1) arc (152.83:180:1.1) -- cycle;
\draw [color=wwttqq] (0,0)-- (3.6,0);
\draw [color=wwttqq] (3.6,0)-- (0.95,1.36);
\draw [color=wwttqq] (0.95,1.36)-- (0,0);
\draw [shift={(0,0)},color=wwttqq] (0:0.8) arc (0:55.06:0.8);
\draw [shift={(0,0)},color=wwttqq] (0:0.7) arc (0:55.06:0.7);
\draw[color=wwttqq] (0.38,0.14) node {$2\beta$};
\draw[color=wwttqq] (2.68,0.16) node {$\beta$};
\end{tikzpicture}}
\\[2pt]
a & b & c \\[5pt]
{\fontsize{9pt}{9pt}\selectfont\begin{tikzpicture}[line cap=round,line join=round,>=triangle 45,x=1.0cm,y=1.0cm]
\clip(-0.4,-0.1) rectangle (2.08,2.52);
\fill[color=ffzzzz,fill=ffzzzz,fill opacity=0.2] (0.25,0) -- (1.9,0) -- (0,2.22) -- cycle;
\draw [shift={(0.25,0)},color=wwttqq,fill=wwttqq,fill opacity=0.2] (0,0) -- (0:0.6) arc (0:96.43:0.6) -- cycle;
\draw [shift={(1.9,0)},color=wwttqq,fill=wwttqq,fill opacity=0.2] (0,0) -- (130.56:0.7) arc (130.56:180:0.7) -- cycle;
\draw [color=wwttqq] (0.25,0)-- (1.9,0);
\draw [color=wwttqq] (1.9,0)-- (0,2.22);
\draw [color=wwttqq] (0,2.22)-- (0.25,0);
\draw [shift={(0.25,0)},color=wwttqq] (0:0.6) arc (0:96.43:0.6);
\draw [shift={(0.25,0)},color=wwttqq] (0:0.5) arc (0:96.43:0.5);
\draw[color=wwttqq] (0.44,0.16) node {$2\beta$};
\draw[color=wwttqq] (1.44,0.16) node {$\beta$};
\end{tikzpicture}}
&
{\fontsize{9pt}{9pt}\selectfont\begin{tikzpicture}[line cap=round,line join=round,>=triangle 45,x=1.0cm,y=1.0cm]
\clip(-0.32,-0.1) rectangle (4.16,0.92);
\fill[color=ffzzzz,fill=ffzzzz,fill opacity=0.2] (0,0) -- (4,0) -- (1.52,0.88) -- cycle;
\draw [shift={(0,0)},color=wwttqq,fill=wwttqq,fill opacity=0.2] (0,0) -- (0:1.2) arc (0:30.07:1.2) -- cycle;
\draw [shift={(4,0)},color=wwttqq,fill=wwttqq,fill opacity=0.2] (0,0) -- (160.46:1.5) arc (160.46:180:1.5) -- cycle;
\draw [shift={(1.52,0.88)},color=wwttqq,fill=wwttqq,fill opacity=0.2] (0,0) -- (-149.93:0.5) arc (-149.93:-19.54:0.5) -- cycle;
\draw [color=wwttqq] (0,0)-- (4,0);
\draw [color=wwttqq] (4,0)-- (1.52,0.88);
\draw [color=wwttqq] (1.52,0.88)-- (0,0);
\draw [shift={(0,0)},color=wwttqq] (0:1.2) arc (0:30.07:1.2);
\draw [shift={(0,0)},color=wwttqq] (0:1.1) arc (0:30.07:1.1);
\draw [shift={(0,0)},color=wwttqq] (0:1) arc (0:30.07:1);
\draw [shift={(4,0)},color=wwttqq] (160.46:1.5) arc (160.46:180:1.5);
\draw [shift={(4,0)},color=wwttqq] (160.46:1.4) arc (160.46:180:1.4);
\draw[color=wwttqq] (0.74,0.16) node {$30^\circ$};
\draw[color=wwttqq] (2.92,0.16) node {$20^\circ$};
\draw[color=wwttqq] (1.62,0.62) node {$130^\circ$};
\end{tikzpicture}}
&
{\fontsize{9pt}{9pt}\selectfont\begin{tikzpicture}[line cap=round,line join=round,>=triangle 45,x=1.0cm,y=1.0cm]
\clip(-0.52,-0.12) rectangle (4.18,2.22);
\fill[color=ffzzzz,fill=ffzzzz,fill opacity=0.2] (0,0) -- (4,0) -- (1.26,1.94) -- cycle;
\draw [shift={(0,0)},color=wwttqq,fill=wwttqq,fill opacity=0.2] (0,0) -- (0:1.2) arc (0:57:1.2) -- cycle;
\draw [shift={(4,0)},color=wwttqq,fill=wwttqq,fill opacity=0.2] (0,0) -- (144.7:1.5) arc (144.7:180:1.5) -- cycle;
\draw [color=wwttqq] (0,0)-- (4,0);
\draw [color=wwttqq] (4,0)-- (1.26,1.94);
\draw [color=wwttqq] (1.26,1.94)-- (0,0);
\draw[color=wwttqq] (0.66,0.3) node {$\frac{n+1}{n}\beta$};
\draw[color=wwttqq] (2.74,0.3) node {$\beta$};
\end{tikzpicture}}
\\[2pt]
d & e & f
\end{tabular}
\caption{Cakes which can be nicely cut into $2$ pieces}
\label{figex}
\end{figure}

\begin{example} \label{ex1} Cakes of indicated shapes (see Figure~\ref{figex}) can be nicely cut into $\mathbf{2}$ pieces:

\noindent
\begin{tabular}{p{2.5cm}p{2.5cm}p{3.5cm}p{3.5cm}}
{\bf (a)} $\alpha=90^\circ$;        & {\bf (b)} $\alpha=3\beta$; &
{\bf (c)} $\alpha=2\beta<90^\circ$; & {\bf (d)} $\alpha=2\beta>90^\circ$; \\
\multicolumn{2}{l}{{\bf (e)} $\alpha=30^\circ$, $\beta=20^\circ$, $\gamma=130^\circ$;} &
\multicolumn{2}{l}{{\bf (f)} $\alpha=\frac{n+1}{n}\beta$, where $n$ is an integer.}
\end{tabular}
\end{example}

We are not going to do the reader out of a pleasure to cut the cakes in Examples~(a)--(d)
nicely himself (herself) by means of straight cuts. Answers are given 
in \S\ref{final}.

Let us present two ways of cutting the cake in Example~(e), which is probably new.


\begin{proof}[\bf 1. A wheel-shaped nice cut]
A cake with $\alpha:\beta=3:2$ can be nicely cut as follows, see Figure~\ref{fig9emore}.
Take a broken line $AKLMNB$ with $5$ equal sides and equal angles
$180^\circ-\alpha+\beta$ between each pair of consecutive sides.
Clearly, the hexagon $AKLMNB$ is mirror-symmetric. Since the sum of the angles of the hexagon is $720^\circ$ it follows that $\angle BAK=\angle ABN=\beta$. Thus $N\in BC$ (assume that both points $N$ and $C$ are above the line $AB$).
The cut $AKLMN$ is nice. Indeed, $\angle CNM=\alpha-\beta=\angle CAK$. Thus the hexagon $AKLMNC$ is also mirror-symmetric and $AC=NC$. Rotating the ''wheel'' $AKLMNB$ one packs the cake into the box.
\end{proof}

\begin{proof}[\bf 2. A gear-shaped nice cut] 
A cake with $\alpha:\beta=3:2$ can be also nicely cut as follows, see Figure~\ref{fig2}.
Take a saw-shaped broken line $AKLMNB$ with $5$ equal sides and angles $\angle AKL=\angle LMN=180^\circ-\beta$,
$\angle KLM=\angle MNB=180^\circ-\alpha$. Then $K\in AB$ and $N\in BC$. The cut $KLMN$ is nice:
rotating the ''gear'' $AKLMNB$ one packs the cake into the box.
\end{proof}

There is even one more way to cut the cake of Example~(e). The reader may wish to find it himself (herself) or see the answer in \S\ref{final}.

Similar wheel- and gear-shaped cuts nicely cut the cake of Example~(f), see Figures~\ref{wheel} and~\ref{saw}.


Conversely, if the cut in Figure~\ref{fig9emore} is nice then $\alpha:\beta=3:2$. Indeed, assume that 
the piece $AKLMNB$ is put into the box by a rotation $R$ such that $R(B)=N=A'$ and $C=C'$. Then $R(N)=M$, $R(M)=L$, $R(L)=K$, $R(K)=A$. A computation shows that the hexagon $AKLMNB$ has $2$ angles of size $\beta$ and $4$ angles of size $180^\circ-\alpha+\beta$.
Since their sum is $720^\circ$ it follows that $2\alpha-3\beta=0$. This is a particular case of the following theorem
(announced in \cite{PSF}):





\begin{theorem}\label{th1} If a cake can be nicely cut into $\mathbf{2}$ pieces then $k\alpha+l\beta+m\gamma=0$ for some integers $k$, $l$ and $m$, not vanishing simultaneously.
\end{theorem}

Thus \emph{almost all} cakes cannot be nicely cut into $\mathbf{2}$ pieces (e.g.,
a cake with $\alpha=\sqrt{2}\,^\circ, \beta=\sqrt{3}\,^\circ, \gamma=180-\sqrt{2}-\sqrt{3}\,^\circ$).

One can try to prove Theorem~\ref{th1} analogously to the above proof that $2\alpha-3\beta=0$ in Figure~\ref{fig9emore}.
But this attempt leads to a huge exhaustion. In this paper we give a short elementary proof of the theorem.

\begin{figure}[phtb]
{\fontsize{9pt}{9pt}\selectfont\begin{tikzpicture}[line cap=round,line join=round,>=triangle 45,x=1.0cm,y=1.0cm]
\clip(0.54,-0.9) rectangle (12.54,3.14);
\fill[color=ffzzzz,fill=ffzzzz,fill opacity=0.2] (0.74,0.3) -- (3.04,2.92) -- (5.31,0.28) -- (4.24,0.89) -- (3.02,1.1) -- (1.8,0.89) -- cycle;
\fill[color=wwttqq,fill=wwttqq,fill opacity=0.4] (0.74,0.3) -- (1.8,0.89) -- (3.02,1.1) -- (4.24,0.89) -- (5.31,0.28) -- (6.12,-0.66) -- cycle;
\fill[color=ffzzzz,fill=ffzzzz,fill opacity=0.2] (12.44,0.3) -- (10.14,2.92) -- (7.86,0.28) -- (8.94,0.89) -- (10.15,1.1) -- (11.37,0.89) -- cycle;
\fill[color=wwttqq,fill=wwttqq,fill opacity=0.4] (12.44,0.3) -- (11.37,0.89) -- (10.15,1.1) -- (8.94,0.89) -- (7.86,0.28) -- (7.06,-0.66) -- cycle;
\draw [shift={(1.8,0.89)},color=wwttqq,fill=wwttqq,fill opacity=0.25] (0,0) -- (-150.79:0.3) arc (-150.79:9.69:0.3) -- cycle;
\draw [shift={(3.02,1.1)},color=wwttqq,fill=wwttqq,fill opacity=0.25] (0,0) -- (-170.31:0.3) arc (-170.31:-9.97:0.3) -- cycle;
\draw [shift={(4.24,0.89)},color=wwttqq,fill=wwttqq,fill opacity=0.25] (0,0) -- (170.03:0.3) arc (170.03:330.37:0.3) -- cycle;
\draw [shift={(5.31,0.28)},color=wwttqq,fill=wwttqq,fill opacity=0.25] (0,0) -- (150.37:0.3) arc (150.37:310.71:0.3) -- cycle;
\draw [shift={(7.86,0.28)},color=wwttqq,fill=wwttqq,fill opacity=0.25] (0,0) -- (-130.71:0.3) arc (-130.71:29.63:0.3) -- cycle;
\draw [shift={(8.94,0.89)},color=wwttqq,fill=wwttqq,fill opacity=0.25] (0,0) -- (-150.37:0.3) arc (-150.37:9.97:0.3) -- cycle;
\draw [shift={(10.15,1.1)},color=wwttqq,fill=wwttqq,fill opacity=0.25] (0,0) -- (-170.03:0.3) arc (-170.03:-9.69:0.3) -- cycle;
\draw [shift={(11.37,0.89)},color=wwttqq,fill=wwttqq,fill opacity=0.25] (0,0) -- (170.31:0.3) arc (170.31:330.79:0.3) -- cycle;
\draw [color=wwttqq] (0.74,0.3)-- (3.04,2.92);
\draw [color=wwttqq] (1.8,1.64) -- (1.93,1.52);
\draw [color=wwttqq] (1.85,1.7) -- (1.98,1.58);
\draw [color=wwttqq] (3.04,2.92)-- (5.31,0.28);
\draw [color=wwttqq] (4.22,1.68) -- (4.09,1.57);
\draw [color=wwttqq] (4.27,1.63) -- (4.13,1.51);
\draw [color=wwttqq] (0.74,0.3)-- (1.8,0.89);
\draw [color=wwttqq] (1.23,0.68) -- (1.32,0.52);
\draw [color=wwttqq] (1.8,0.89)-- (3.02,1.1);
\draw [color=wwttqq] (2.4,1.09) -- (2.43,0.91);
\draw [color=wwttqq] (3.02,1.1)-- (4.24,0.89);
\draw [color=wwttqq] (3.65,1.08) -- (3.61,0.91);
\draw [color=wwttqq] (4.24,0.89)-- (5.31,0.28);
\draw [color=wwttqq] (4.82,0.66) -- (4.73,0.5);
\draw [color=wwttqq] (5.31,0.28)-- (6.12,-0.66);
\draw [color=wwttqq] (5.79,-0.13) -- (5.65,-0.25);
\draw [color=wwttqq] (6.12,-0.66)-- (0.74,0.3);
\draw [color=wwttqq] (12.44,0.3)-- (10.14,2.92);
\draw [color=wwttqq] (11.24,1.52) -- (11.38,1.64);
\draw [color=wwttqq] (11.19,1.58) -- (11.33,1.7);
\draw [color=wwttqq] (10.14,2.92)-- (7.86,0.28);
\draw [color=wwttqq] (9.09,1.57) -- (8.95,1.68);
\draw [color=wwttqq] (9.04,1.51) -- (8.91,1.63);
\draw [color=wwttqq] (12.44,0.3)-- (11.37,0.89);
\draw [color=wwttqq] (11.86,0.52) -- (11.95,0.68);
\draw [color=wwttqq] (11.37,0.89)-- (10.15,1.1);
\draw [color=wwttqq] (10.75,0.91) -- (10.78,1.09);
\draw [color=wwttqq] (10.15,1.1)-- (8.94,0.89);
\draw [color=wwttqq] (9.56,0.91) -- (9.53,1.08);
\draw [color=wwttqq] (8.94,0.89)-- (7.86,0.28);
\draw [color=wwttqq] (8.44,0.5) -- (8.35,0.66);
\draw [color=wwttqq] (7.86,0.28)-- (7.06,-0.66);
\draw [color=wwttqq] (7.53,-0.25) -- (7.39,-0.13);
\draw [color=wwttqq] (7.06,-0.66)-- (12.44,0.3);
\draw [shift={(3.01,-2.52)},color=wwttqq]  plot[domain=1.19:1.52,variable=\t]({1*3.05*cos(\t r)+0*3.05*sin(\t r)},{0*3.05*cos(\t r)+1*3.05*sin(\t r)});
\draw [->,color=wwttqq] (3.18,0.53) -- (2.91,0.53);
\fill [color=wwttqq] (0.74,0.3) circle (1.0pt);
\draw[color=wwttqq] (0.66,0.58) node {$A$};
\fill [color=wwttqq] (6.12,-0.66) circle (1.0pt);
\draw[color=wwttqq] (6.26,-0.41) node {$B$};
\fill [color=wwttqq] (3.04,2.92) circle (1.0pt);
\draw[color=wwttqq] (2.65,2.82) node {$C$};
\fill [color=wwttqq] (5.31,0.28) circle (1.0pt);
\draw[color=wwttqq] (5.46,0.47) node {$N$};
\fill [color=wwttqq] (3.02,1.1) circle (1.0pt);
\draw[color=wwttqq] (2.99,1.36) node {$L$};
\fill [color=wwttqq] (4.24,0.89) circle (1.0pt);
\draw[color=wwttqq] (4.21,1.16) node {$M$};
\fill [color=wwttqq] (1.8,0.89) circle (1.0pt);
\draw[color=wwttqq] (1.84,1.14) node {$K$};
\fill [color=wwttqq] (10.14,2.92) circle (1.0pt);
\fill [color=wwttqq] (12.44,0.3) circle (1.0pt);
\fill [color=wwttqq] (11.37,0.89) circle (1.0pt);
\fill [color=wwttqq] (10.15,1.1) circle (1.0pt);
\fill [color=wwttqq] (8.94,0.89) circle (1.0pt);
\fill [color=wwttqq] (7.86,0.28) circle (1.0pt);
\fill [color=wwttqq] (7.06,-0.66) circle (1.0pt);
\end{tikzpicture}}
\caption{A wheel-shaped nice cut}
\label{fig9emore}
\end{figure}
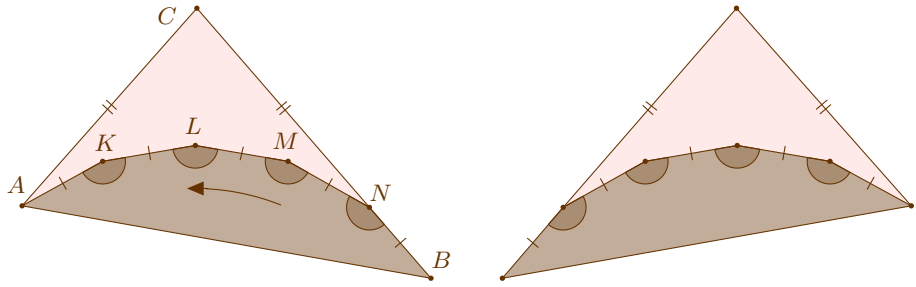

\begin{figure}[phtb]
{\fontsize{9pt}{9pt}\selectfont\begin{tikzpicture}[line cap=round,line join=round,>=triangle 45,x=1.0cm,y=1.0cm]
\clip(0.54,-0.84) rectangle (12.56,3.04);
\fill[color=ffzzzz,fill=ffzzzz,fill opacity=0.2] (0.74,0.3) -- (1.9573,0.084) -- (3.0288,0.6812) -- (4.0944,0.0734) -- (5.3137,0.2772) -- (3.04,2.92) -- cycle;
\fill[color=wwttqq,fill=wwttqq,fill opacity=0.4] (1.9573,0.084) -- (3.0288,0.6812) -- (4.0944,0.0734) -- (5.3137,0.2772) -- (6.12,-0.66) -- cycle;
\draw [shift={(3.0288,0.6812)},color=wwttqq,fill=wwttqq,fill opacity=0.25] (0,0) -- (-150.8685:0.3) arc (-150.8685:-29.7013:0.3) -- cycle;
\draw [shift={(5.3137,0.2772)},color=wwttqq,fill=wwttqq,fill opacity=0.25] (0,0) -- (-170.5133:0.3) arc (-170.5133:-49.2934:0.3) -- cycle;
\draw [shift={(4.0944,0.0734)},color=ffzzzz,fill=ffzzzz,fill opacity=0.25] (0,0) -- (9.4867:0.4) arc (9.4867:150.2987:0.4) -- cycle;
\draw [shift={(1.9573,0.084)},color=ffzzzz,fill=ffzzzz,fill opacity=0.25] (0,0) -- (29.1315:0.4) arc (29.1315:169.9393:0.4) -- cycle;
\fill[color=ffzzzz,fill=ffzzzz,fill opacity=0.2] (12.4354,0.3) -- (11.2181,0.084) -- (10.1466,0.6812) -- (9.081,0.0734) -- (7.8617,0.2772) -- (10.1354,2.92) -- cycle;
\fill[color=wwttqq,fill=wwttqq,fill opacity=0.4] (11.2181,0.084) -- (10.1466,0.6812) -- (9.081,0.0734) -- (7.8617,0.2772) -- (7.0554,-0.66) -- cycle;
\draw [shift={(11.2181,0.084)},color=ffzzzz,fill=ffzzzz,fill opacity=0.25] (0,0) -- (10.0607:0.4) arc (10.0607:150.8685:0.4) -- cycle;
\draw [shift={(9.081,0.0734)},color=ffzzzz,fill=ffzzzz,fill opacity=0.25] (0,0) -- (29.7013:0.4) arc (29.7013:170.5133:0.4) -- cycle;
\draw [shift={(7.8617,0.2772)},color=wwttqq,fill=wwttqq,fill opacity=0.25] (0,0) -- (-130.7066:0.3) arc (-130.7066:-9.4867:0.3) -- cycle;
\draw [shift={(10.1466,0.6812)},color=wwttqq,fill=wwttqq,fill opacity=0.25] (0,0) -- (-150.2987:0.3) arc (-150.2987:-29.1315:0.3) -- cycle;
\draw [color=wwttqq] (0.74,0.3)-- (1.9573,0.084);
\draw [color=wwttqq] (1.3644,0.2806) -- (1.3329,0.1034);
\draw [color=wwttqq] (1.9573,0.084)-- (3.0288,0.6812);
\draw [color=wwttqq] (2.4492,0.4612) -- (2.5369,0.304);
\draw [color=wwttqq] (3.0288,0.6812)-- (4.0944,0.0734);
\draw [color=wwttqq] (3.6062,0.4555) -- (3.517,0.2991);
\draw [color=wwttqq] (4.0944,0.0734)-- (5.3137,0.2772);
\draw [color=wwttqq] (4.6892,0.2641) -- (4.7189,0.0865);
\draw [color=wwttqq] (5.3137,0.2772)-- (3.04,2.92);
\draw [color=wwttqq] (4.1315,1.5134) -- (4.2679,1.6307);
\draw [color=wwttqq] (4.0858,1.5664) -- (4.2223,1.6838);
\draw [color=wwttqq] (3.04,2.92)-- (0.74,0.3);
\draw [color=wwttqq] (1.9807,1.5769) -- (1.8455,1.6957);
\draw [color=wwttqq] (1.9345,1.5243) -- (1.7993,1.6431);
\draw [color=wwttqq] (1.9573,0.084)-- (3.0288,0.6812);
\draw [color=wwttqq] (3.0288,0.6812)-- (4.0944,0.0734);
\draw [color=wwttqq] (4.0944,0.0734)-- (5.3137,0.2772);
\draw [color=wwttqq] (5.3137,0.2772)-- (6.12,-0.66);
\draw [color=wwttqq] (5.7851,-0.1327) -- (5.6486,-0.2501);
\draw [color=wwttqq] (6.12,-0.66)-- (1.9573,0.084);
\draw [shift={(4.0944,0.0734)},color=ffzzzz] (9.4867:0.4) arc (9.4867:150.2987:0.4);
\draw [shift={(4.0944,0.0734)},color=ffzzzz] (9.4867:0.3) arc (9.4867:150.2987:0.3);
\draw [shift={(1.9573,0.084)},color=ffzzzz] (29.1315:0.4) arc (29.1315:169.9393:0.4);
\draw [shift={(1.9573,0.084)},color=ffzzzz] (29.1315:0.3) arc (29.1315:169.9393:0.3);
\draw [color=wwttqq] (12.4354,0.3)-- (11.2181,0.084);
\draw [color=wwttqq] (11.8425,0.1034) -- (11.811,0.2806);
\draw [color=wwttqq] (11.2181,0.084)-- (10.1466,0.6812);
\draw [color=wwttqq] (10.6385,0.304) -- (10.7262,0.4612);
\draw [color=wwttqq] (10.1466,0.6812)-- (9.081,0.0734);
\draw [color=wwttqq] (9.6584,0.2991) -- (9.5692,0.4555);
\draw [color=wwttqq] (9.081,0.0734)-- (7.8617,0.2772);
\draw [color=wwttqq] (8.4565,0.0865) -- (8.4862,0.2641);
\draw [color=wwttqq] (7.8617,0.2772)-- (10.1354,2.92);
\draw [color=wwttqq] (8.9075,1.6307) -- (9.0439,1.5134);
\draw [color=wwttqq] (8.9531,1.6838) -- (9.0896,1.5664);
\draw [color=wwttqq] (10.1354,2.92)-- (12.4354,0.3);
\draw [color=wwttqq] (11.3299,1.6957) -- (11.1947,1.5769);
\draw [color=wwttqq] (11.3761,1.6431) -- (11.2409,1.5243);
\draw [color=wwttqq] (11.2181,0.084)-- (10.1466,0.6812);
\draw [color=wwttqq] (10.1466,0.6812)-- (9.081,0.0734);
\draw [color=wwttqq] (9.081,0.0734)-- (7.8617,0.2772);
\draw [color=wwttqq] (7.8617,0.2772)-- (7.0554,-0.66);
\draw [color=wwttqq] (7.5268,-0.2501) -- (7.3903,-0.1327);
\draw [color=wwttqq] (7.0554,-0.66)-- (11.2181,0.084);
\draw [shift={(2.9948,-6.1275)},color=wwttqq]  plot[domain=1.2234:1.5658,variable=\t]({1*5.8048*cos(\t r)+0*5.8048*sin(\t r)},{0*5.8048*cos(\t r)+1*5.8048*sin(\t r)});
\draw [->,color=wwttqq] (3.28,-0.34) -- (3.0237,-0.3401);
\draw [shift={(11.2181,0.084)},color=ffzzzz] (10.0607:0.4) arc (10.0607:150.8685:0.4);
\draw [shift={(11.2181,0.084)},color=ffzzzz] (10.0607:0.3) arc (10.0607:150.8685:0.3);
\draw [shift={(9.081,0.0734)},color=ffzzzz] (29.7013:0.4) arc (29.7013:170.5133:0.4);
\draw [shift={(9.081,0.0734)},color=ffzzzz] (29.7013:0.3) arc (29.7013:170.5133:0.3);
\fill [color=wwttqq] (0.74,0.3) circle (1.0pt);
\draw[color=wwttqq] (0.66,0.54) node {$A$};
\fill [color=wwttqq] (6.12,-0.66) circle (1.0pt);
\draw[color=wwttqq] (6.28,-0.42) node {$B$};
\fill [color=wwttqq] (3.04,2.92) circle (1.0pt);
\draw[color=wwttqq] (2.72,2.92) node {$C$};
\fill [color=wwttqq] (5.3137,0.2772) circle (1.0pt);
\draw[color=wwttqq] (5.48,0.52) node {$N$};
\fill [color=wwttqq] (1.9573,0.084) circle (1.0pt);
\draw[color=wwttqq] (1.78,-0.14) node {$K$};
\fill [color=wwttqq] (3.0288,0.6812) circle (1.0pt);
\draw[color=wwttqq] (3.16,0.92) node {$L$};
\fill [color=wwttqq] (4.0944,0.0734) circle (1.0pt);
\draw[color=wwttqq] (4.06,-0.12) node {$M$};
\fill [color=wwttqq] (12.4354,0.3) circle (1.0pt);
\fill [color=wwttqq] (11.2181,0.084) circle (1.0pt);
\fill [color=wwttqq] (10.1466,0.6812) circle (1.0pt);
\fill [color=wwttqq] (9.081,0.0734) circle (1.0pt);
\fill [color=wwttqq] (7.8617,0.2772) circle (1.0pt);
\fill [color=wwttqq] (10.1354,2.92) circle (1.0pt);
\fill [color=wwttqq] (7.0554,-0.66) circle (1.0pt);
\end{tikzpicture}}
\caption{A gear-shaped nice cut}
\label{fig2}
\end{figure}
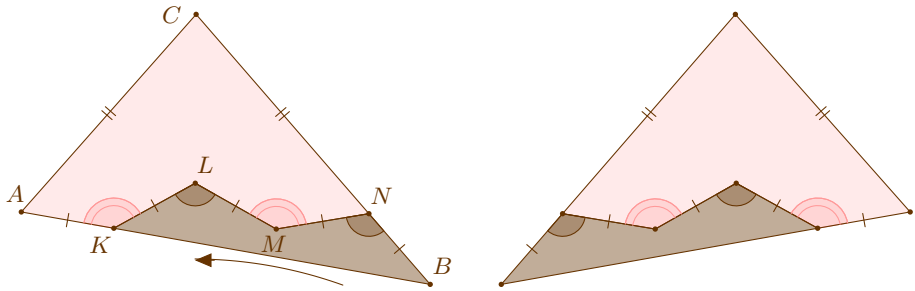

\begin{figure}[phtb]
{\fontsize{9pt}{9pt}\selectfont\begin{tikzpicture}[line cap=round,line join=round,>=triangle 45,x=1.0cm,y=1.0cm]
\clip(-1.98,1.5) rectangle (7.52,5.84);
\fill[color=wwttqq,fill=wwttqq,fill opacity=0.4] (-1.39106,3.23815) -- (-1.02945,3.59465) -- (-0.59054,3.85003) -- (-0.10191,3.98823) -- (0.40574,4.00057) -- (0.9005,3.88628) -- (1.3513,3.65254) -- (1.72981,3.31403) -- (2.01224,2.89202) -- cycle;
\fill[color=ffzzzz,fill=ffzzzz,fill opacity=0.2] (-1.39106,3.23815) -- (0.1096,5.73493) -- (1.72981,3.31403) -- (1.3513,3.65254) -- (0.9005,3.88628) -- (0.40574,4.00057) -- (-0.10191,3.98823) -- (-0.59054,3.85003) -- (-1.02945,3.59465) -- cycle;
\fill[color=ffzzzz,fill=ffzzzz,fill opacity=0.2] (6.65356,3.23815) -- (5.1529,5.73493) -- (3.53269,3.31403) -- (3.9112,3.65254) -- (4.362,3.88628) -- (4.85676,4.00057) -- (5.36441,3.98823) -- (5.85304,3.85003) -- (6.29195,3.59465) -- cycle;
\fill[color=wwttqq,fill=wwttqq,fill opacity=0.4] (6.65356,3.23815) -- (6.29195,3.59465) -- (5.85304,3.85003) -- (5.36441,3.98823) -- (4.85676,4.00057) -- (4.362,3.88628) -- (3.9112,3.65254) -- (3.53269,3.31403) -- (3.25026,2.89202) -- cycle;
\draw [color=wwttqq] (-1.39106,3.23815)-- (-1.02945,3.59465);
\draw [color=wwttqq] (-1.02945,3.59465)-- (-0.59054,3.85003);
\draw [color=wwttqq] (-0.59054,3.85003)-- (-0.10191,3.98823);
\draw [color=wwttqq] (-0.10191,3.98823)-- (0.40574,4.00057);
\draw [color=wwttqq] (0.40574,4.00057)-- (0.9005,3.88628);
\draw [color=wwttqq] (0.9005,3.88628)-- (1.3513,3.65254);
\draw [color=wwttqq] (1.3513,3.65254)-- (1.72981,3.31403);
\draw [color=wwttqq] (1.72981,3.31403)-- (2.01224,2.89202);
\draw [color=wwttqq] (2.01224,2.89202)-- (-1.39106,3.23815);
\draw [color=wwttqq] (-1.39106,3.23815)-- (0.1096,5.73493);
\draw [color=wwttqq] (0.1096,5.73493)-- (1.72981,3.31403);
\draw [dotted,color=wwttqq] (0.20076,1.98519) circle (0.29236cm);
\draw [dotted,color=wwttqq] (-1.39106,3.23815)-- (-1.65264,2.80292);
\draw [dotted,color=wwttqq] (-1.65264,2.80292)-- (-1.79778,2.3163);
\draw [dotted,color=wwttqq] (2.01224,2.89202)-- (2.18085,2.41303);
\draw [color=wwttqq] (6.65356,3.23815)-- (5.1529,5.73493);
\draw [color=wwttqq] (5.1529,5.73493)-- (3.53269,3.31403);
\draw [color=wwttqq] (6.65356,3.23815)-- (6.29195,3.59465);
\draw [color=wwttqq] (6.29195,3.59465)-- (5.85304,3.85003);
\draw [color=wwttqq] (5.85304,3.85003)-- (5.36441,3.98823);
\draw [color=wwttqq] (5.36441,3.98823)-- (4.85676,4.00057);
\draw [color=wwttqq] (4.85676,4.00057)-- (4.362,3.88628);
\draw [color=wwttqq] (4.362,3.88628)-- (3.9112,3.65254);
\draw [color=wwttqq] (3.9112,3.65254)-- (3.53269,3.31403);
\draw [color=wwttqq] (3.53269,3.31403)-- (3.25026,2.89202);
\draw [color=wwttqq] (3.25026,2.89202)-- (6.65356,3.23815);
\draw [dotted,color=wwttqq] (6.65356,3.23815)-- (6.91514,2.80292);
\draw [dotted,color=wwttqq] (6.91514,2.80292)-- (7.06028,2.3163);
\draw [dotted,color=wwttqq] (3.25026,2.89202)-- (3.08165,2.41303);
\draw [dotted,color=wwttqq] (5.06174,1.98519) circle (0.29236cm);
\draw [shift={(0.20076,1.98519)},color=wwttqq]  plot[domain=0.34014:2.66392,variable=\t]({1*0.61821*cos(\t r)+0*0.61821*sin(\t r)},{0*0.61821*cos(\t r)+1*0.61821*sin(\t r)});
\draw [->,color=wwttqq] (-0.24875,2.38) -- (-0.34824,2.26939);
\fill [color=wwttqq] (2.18085,2.41303) circle (1.0pt);
\fill [color=wwttqq] (2.01224,2.89202) circle (1.0pt);
\fill [color=wwttqq] (1.72981,3.31403) circle (1.0pt);
\fill [color=wwttqq] (1.3513,3.65254) circle (1.0pt);
\fill [color=wwttqq] (0.9005,3.88628) circle (1.0pt);
\fill [color=wwttqq] (0.40574,4.00057) circle (1.0pt);
\fill [color=wwttqq] (-0.10191,3.98823) circle (1.0pt);
\fill [color=wwttqq] (-0.59054,3.85003) circle (1.0pt);
\fill [color=wwttqq] (-1.02945,3.59465) circle (1.0pt);
\fill [color=wwttqq] (-1.39106,3.23815) circle (1.0pt);
\fill [color=wwttqq] (-1.65264,2.80292) circle (1.0pt);
\fill [color=wwttqq] (-1.79778,2.3163) circle (1.0pt);
\fill [color=wwttqq] (6.29195,3.59465) circle (1.0pt);
\fill [color=wwttqq] (5.85304,3.85003) circle (1.0pt);
\fill [color=wwttqq] (5.36441,3.98823) circle (1.0pt);
\fill [color=wwttqq] (4.85676,4.00057) circle (1.0pt);
\fill [color=wwttqq] (4.362,3.88628) circle (1.0pt);
\fill [color=wwttqq] (3.9112,3.65254) circle (1.0pt);
\fill [color=wwttqq] (3.53269,3.31403) circle (1.0pt);
\fill [color=wwttqq] (6.65356,3.23815) circle (1.0pt);
\fill [color=wwttqq] (6.91514,2.80292) circle (1.0pt);
\fill [color=wwttqq] (7.06028,2.3163) circle (1.0pt);
\fill [color=wwttqq] (3.25026,2.89202) circle (1.0pt);
\fill [color=wwttqq] (3.08165,2.41303) circle (1.0pt);
\end{tikzpicture}}
\caption{A wheel-shaped nice cut again}
\label{wheel}
\end{figure}
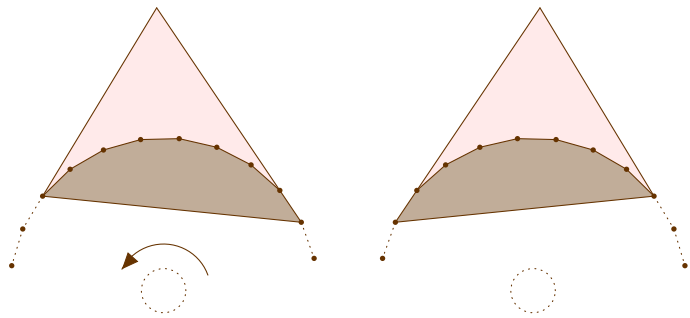

\begin{figure}[phtb]
{\fontsize{9pt}{9pt}\selectfont\begin{tikzpicture}[scale=1.3,line cap=round,line join=round,>=triangle 45,x=1.0cm,y=1.0cm]
\clip(0.08,1.5) rectangle (9.06,5.84);
\fill[color=wwttqq,fill=wwttqq,fill opacity=0.4] (1.52,3.83) -- (1.72,4.05) -- (2,3.96) -- (2.25,4.12) -- (2.5,3.96) -- (2.78,4.05) -- (2.98,3.83) -- (3.28,3.85) -- (3.42,3.59) -- cycle;
\fill[color=ffzzzz,fill=ffzzzz,fill opacity=0.2] (3.28,3.85) -- (2.25,5.72) -- (1.22,3.85) -- (1.52,3.83) -- (1.72,4.05) -- (2,3.96) -- (2.25,4.12) -- (2.5,3.96) -- (2.78,4.05) -- (2.98,3.83) -- cycle;
\fill[color=ffzzzz,fill=ffzzzz,fill opacity=0.2] (5.84,3.85) -- (6.87,5.72) -- (7.9,3.85) -- (7.6,3.83) -- (7.4,4.05) -- (7.12,3.96) -- (6.87,4.12) -- (6.62,3.96) -- (6.34,4.05) -- (6.14,3.83) -- cycle;
\fill[color=wwttqq,fill=wwttqq,fill opacity=0.4] (7.6,3.83) -- (7.4,4.05) -- (7.12,3.96) -- (6.87,4.12) -- (6.62,3.96) -- (6.34,4.05) -- (6.14,3.83) -- (5.84,3.85) -- (5.7,3.59) -- cycle;
\draw [color=wwttqq] (1.52,3.83)-- (1.72,4.05);
\draw [color=wwttqq] (1.72,4.05)-- (2,3.96);
\draw [color=wwttqq] (2,3.96)-- (2.25,4.12);
\draw [color=wwttqq] (2.25,4.12)-- (2.5,3.96);
\draw [color=wwttqq] (2.5,3.96)-- (2.78,4.05);
\draw [color=wwttqq] (2.78,4.05)-- (2.98,3.83);
\draw [color=wwttqq] (2.98,3.83)-- (3.28,3.85);
\draw [color=wwttqq] (3.28,3.85)-- (3.42,3.59);
\draw [color=wwttqq] (3.42,3.59)-- (1.52,3.83);
\draw [color=wwttqq] (3.28,3.85)-- (2.25,5.72);
\draw [color=wwttqq] (2.25,5.72)-- (1.22,3.85);
\draw [color=wwttqq] (1.22,3.85)-- (1.52,3.83);
\draw [dotted,color=wwttqq] (0.22,2.64)-- (0.45,2.83);
\draw [dotted,color=wwttqq] (0.45,3.12)-- (0.45,2.83);
\draw [dotted,color=wwttqq] (0.45,3.12)-- (0.71,3.25);
\draw [dotted,color=wwttqq] (0.71,3.25)-- (0.79,3.54);
\draw [dotted,color=wwttqq] (0.79,3.54)-- (1.08,3.59);
\draw [dotted,color=wwttqq] (1.08,3.59)-- (1.22,3.85);
\draw [dotted,color=wwttqq] (3.42,3.59)-- (3.71,3.54);
\draw [dotted,color=wwttqq] (3.71,3.54)-- (3.79,3.25);
\draw [dotted,color=wwttqq] (3.79,3.25)-- (4.05,3.12);
\draw [dotted,color=wwttqq] (4.05,3.12)-- (4.05,2.83);
\draw [dotted,color=wwttqq] (4.05,2.83)-- (4.28,2.64);
\draw [dotted,color=wwttqq] (4.28,2.64)-- (4.21,2.35);
\draw [dotted,color=wwttqq] (4.21,2.35)-- (4.38,2.11);
\draw [dotted,color=wwttqq] (8.04,3.59)-- (7.9,3.85);
\draw [dotted,color=wwttqq] (8.33,3.54)-- (8.04,3.59);
\draw [dotted,color=wwttqq] (8.41,3.25)-- (8.33,3.54);
\draw [dotted,color=wwttqq] (8.67,3.12)-- (8.41,3.25);
\draw [dotted,color=wwttqq] (8.67,3.12)-- (8.67,2.83);
\draw [dotted,color=wwttqq] (8.9,2.64)-- (8.67,2.83);
\draw [dotted,color=wwttqq] (5.7,3.59)-- (5.41,3.54);
\draw [dotted,color=wwttqq] (5.41,3.54)-- (5.33,3.25);
\draw [dotted,color=wwttqq] (5.33,3.25)-- (5.07,3.12);
\draw [dotted,color=wwttqq] (5.07,3.12)-- (5.07,2.83);
\draw [dotted,color=wwttqq] (5.07,2.83)-- (4.84,2.64);
\draw [dotted,color=wwttqq] (4.84,2.64)-- (4.91,2.35);
\draw [dotted,color=wwttqq] (4.91,2.35)-- (4.74,2.11);
\draw [dotted,color=wwttqq] (2.25,1.98) circle (0.32cm);
\draw [dotted,color=wwttqq] (0.29,2.35)-- (0.22,2.64);
\draw [dotted,color=wwttqq] (0.29,2.35)-- (0.12,2.11);
\draw [color=wwttqq] (5.84,3.85)-- (6.87,5.72);
\draw [color=wwttqq] (6.87,5.72)-- (7.9,3.85);
\draw [color=wwttqq] (7.9,3.85)-- (7.6,3.83);
\draw [color=wwttqq] (7.6,3.83)-- (7.4,4.05);
\draw [color=wwttqq] (7.4,4.05)-- (7.12,3.96);
\draw [color=wwttqq] (7.12,3.96)-- (6.87,4.12);
\draw [color=wwttqq] (6.87,4.12)-- (6.62,3.96);
\draw [color=wwttqq] (6.62,3.96)-- (6.34,4.05);
\draw [color=wwttqq] (6.34,4.05)-- (6.14,3.83);
\draw [color=wwttqq] (6.14,3.83)-- (5.84,3.85);
\draw [color=wwttqq] (7.6,3.83)-- (7.4,4.05);
\draw [color=wwttqq] (7.4,4.05)-- (7.12,3.96);
\draw [color=wwttqq] (7.12,3.96)-- (6.87,4.12);
\draw [color=wwttqq] (6.87,4.12)-- (6.62,3.96);
\draw [color=wwttqq] (6.62,3.96)-- (6.34,4.05);
\draw [color=wwttqq] (6.34,4.05)-- (6.14,3.83);
\draw [color=wwttqq] (6.14,3.83)-- (5.84,3.85);
\draw [color=wwttqq] (5.84,3.85)-- (5.7,3.59);
\draw [color=wwttqq] (5.7,3.59)-- (7.6,3.83);
\draw [dotted,color=wwttqq] (8.83,2.35)-- (8.9,2.64);
\draw [dotted,color=wwttqq] (8.83,2.35)-- (9,2.11);
\draw [dotted,color=wwttqq] (6.87,1.98) circle (0.32cm);
\draw [shift={(2.25,1.98)},color=wwttqq]  plot[domain=0.4:2.1,variable=\t]({1*0.6*cos(\t r)+0*0.6*sin(\t r)},{0*0.6*cos(\t r)+1*0.6*sin(\t r)});
\draw [->,color=wwttqq] (1.95,2.5) -- (1.76,2.42);
\end{tikzpicture}}
\caption{A gear-shaped nice cut again}
\label{saw}
\end{figure}
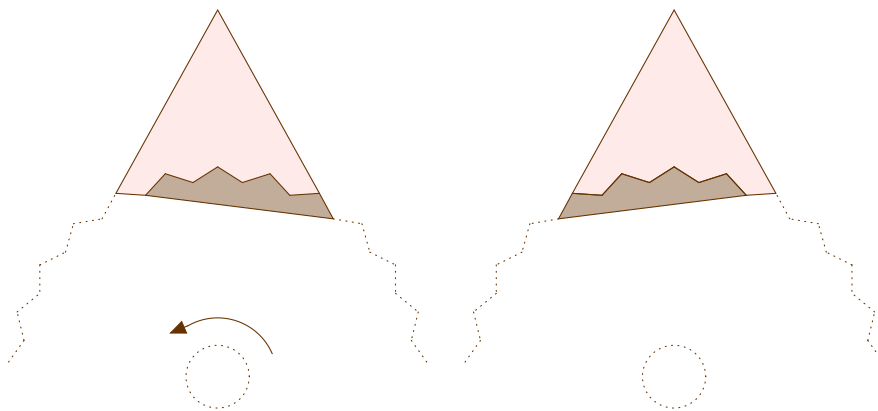

\section{Invariants}

Let us discuss historical background of Theorem~\ref{th1}. Our history begins with

\begin{theorem*}[Bolyai--Gerwien, 1832] Suppose that two planar polygons have the same area. 
Then the first polygon can be cut into finitely many polygonal pieces which can be 
reassembled as the second polygon \textup{(}possibly turning over the pieces\textup{)}.
\end{theorem*}


V.G. Boltyansky comes to his question while proving that the Bolyai-Gerwien theorem remains true if turning over the pieces is not permitted.
Let us sketch the proof. Assume that a piece $P$ is turned over during the assembling. Triangulate the piece. Cut each 
of the obtained triangles into $3$ subpieces as in Example~\ref{ex0}. These subpieces can be reassembled as the mirror image of~$P$.
Thus one gets rid of turnings.

H. Hadwiger and P. Glur go even further. 
They prove that Bolyai--Gerwien theorem remains true if only parallel translations and central symmetries of the pieces are permitted \cite[\S9]{Bol77}. The latter result is optimal: this collection of isometries is in some sense minimal \cite[\S10]{Bol77}.
The proof of optimality uses an \emph{additive invariant}, which
is also the main tool of the paper.



\begin{figure}[htb]
{\fontsize{9pt}{9pt}\selectfont\begin{tikzpicture}[line cap=round,line join=round,>=triangle 45,x=0.66cm,y=0.66cm]
\clip(-2.6,-2.58) rectangle (6.4,4.72);
\fill[line width=0pt,color=ffzzzz,fill=ffzzzz,fill opacity=0.2] (-2.28,2) -- (-1,2) -- (0,1) -- (1,3) -- (4,3) -- (5,2) -- (6,2) -- (6,0) -- (4,-2) -- (1,-2) -- (0,-1) -- (-1.2,1) -- cycle;
\draw [->,line width=1.2pt,color=zzqqcc] (5,4) -- (-1,4);
\draw [->,line width=1.2pt,color=zzqqcc] (-1,2) -- (-2.28,2);
\draw [->,color=wwttqq] (-2.28,2) -- (-1.2,1);
\draw [->,color=wwttqq] (-1.2,1) -- (0,-1);
\draw [->,color=wwttqq] (0,-1) -- (1,-2);
\draw [->,line width=1.2pt,color=zzqqcc] (1,-2) -- (4,-2);
\draw [->,color=wwttqq] (4,-2) -- (6,0);
\draw [->,color=wwttqq] (6,0) -- (6,2);
\draw [->,line width=1.2pt,color=zzqqcc] (6,2) -- (5,2);
\draw [->,color=wwttqq] (5,2) -- (4,3);
\draw [->,line width=1.2pt,color=zzqqcc] (4,3) -- (1,3);
\draw [->,color=wwttqq] (1,3) -- (0,1);
\draw [->,color=wwttqq] (0,1) -- (-1,2);
\draw[color=wwttqq] (-2.46,2.3) node {$G$};
\draw[color=wwttqq] (-0.8,2.28) node {$F$};
\draw[color=wwttqq] (0.36,1.08) node {$E$};
\draw[color=wwttqq] (1.16,3.32) node {$D$};
\draw[color=wwttqq] (4.16,3.24) node {$C$};
\draw[color=wwttqq] (5.04,2.32) node {$B$};
\draw[color=wwttqq] (6.18,2.28) node {$A$};
\draw[color=wwttqq] (6.24,0.22) node {$L$};
\draw[color=wwttqq] (3.84,-1.6) node {$K$};
\draw[color=wwttqq] (1.18,-1.62) node {$J$};
\draw[color=wwttqq] (0.2,-0.7) node {$I$};
\draw[color=wwttqq] (-0.82,1.02) node {$H$};
\draw[color=ffzzzz] (2.44,0.6) node {$P$};
\draw[color=zzqqcc] (2,4.3) node {$f$};
\draw[color=zzqqcc] (-1.62,2.42) node {+};
\draw[color=zzqqcc] (2.62,-2.22) node {$-$};
\draw[color=zzqqcc] (5.58,2.36) node {+};
\draw[color=zzqqcc] (2.52,3.42) node {+};
\end{tikzpicture}}
\caption{Definition of the additive invariant: $J_f(P)=AB+CD+FG-JK$}
\label{figinv}
\end{figure}
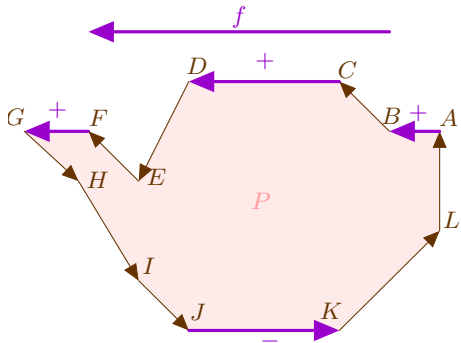

\begin{definition*}[of the additive invariant] \cite[\S10]{Bol77} Let 
$P$ be a polygon in the plane, see Figure~\ref{figinv}.
Orient all the sides of the polygon 
$P$ counterclockwise.
Let $f$ be a {\it directed line}, i.e., a line with a direction on it marked by an arrow.
The {\it additive invariant} $J_f(P)$ is an algebraic sum 
of all the sides of the polygon $P$
parallel to $f$. In this sum the sides having the {\it same}
direction as $f$ (sides  $AB$, $CD$ and $FG$ in Figure~\ref{figinv}), are taken
with positive 
sign, and the ones having the {\it opposite} direction (side
$JK$ in Figure~\ref{figinv}) are taken with negative 
sign. If the polygon
$P$ has no sides parallel to $f$ then set $J_f(P)=0$.
\end{definition*}

Let us establish some basic properties of the additive invariant.

\begin{claim}\label{cl1} \cite[\S10]{Bol77}
\textup{\bf (a)} Suppose that a polygon 
$P$ is cut into several polygons
$P_1,\dots,
P_n$.
Then $J_f(P)=
J_f(P_1)+\dots+
J_f(P_n)$.

\textup{\bf (b)} Suppose that a cake 
$\begin{tikzpicture}[line cap=round,line join=round,>=triangle 45,x=0.11452372190810167cm,y=0.11452372190810167cm]
\clip(-1.75,3.88) rectangle (1.74,5.89);
\draw (-1.72,5.78)-- (0.72,4.38);
\draw (0.72,4.38)-- (1.56,5.78);
\draw (-1.72,5.78)-- (1.56,5.78);
\draw (-1.72,5.78)-- (-1.7,5.46);
\draw (-1.7,5.46)-- (0.72,4.08);
\draw (0.72,4.38)-- (0.72,4.08);
\draw (0.72,4.08)-- (1.56,5.46);
\draw (1.56,5.78)-- (1.56,5.46);
\draw [rotate around={-17.53:(0.13,5.38)}] (0.13,5.38) ellipse (0.03cm and 0.03cm);
\draw (0.42,4.8)-- (0.74,4.8);
\draw (0.74,4.8)-- (0.58,5.04);
\draw (0.58,5.04)-- (0.96,5.1);
\draw (0.96,5.1)-- (0.78,5.32);
\draw (0.78,5.32)-- (1.14,5.42);
\draw (1.14,5.42)-- (0.92,5.56);
\draw (0.92,5.56)-- (1.26,5.64);
\end{tikzpicture}$ is cut into several polygonal pieces which are put into the box 
$\begin{tikzpicture}[line cap=round,line join=round,>=triangle 45,x=0.11711167869707381cm,y=0.11711167869707381cm]
\clip(2.62,3.92) rectangle (6.04,6.09);
\draw (2.66,5.92)-- (5.92,5.92);
\draw (2.66,5.92)-- (3.48,4.52);
\draw (5.92,5.92)-- (3.48,4.52);
\draw (2.66,5.92)-- (2.66,5.44);
\draw (3.48,4.52)-- (3.5,4.02);
\draw (5.92,5.92)-- (5.94,5.46);
\draw (2.66,5.44)-- (3.5,4.02);
\draw (3.5,4.02)-- (5.94,5.46);
\draw (2.95,5.43)-- (5.09,5.44);
\end{tikzpicture}$ all together by means of parallel translations.
Then
$J_f(\begin{tikzpicture}[line cap=round,line join=round,>=triangle 45,x=0.11452372190810167cm,y=0.11452372190810167cm]
\clip(-1.75,3.88) rectangle (1.74,5.89);
\draw (-1.72,5.78)-- (0.72,4.38);
\draw (0.72,4.38)-- (1.56,5.78);
\draw (-1.72,5.78)-- (1.56,5.78);
\draw (-1.72,5.78)-- (-1.7,5.46);
\draw (-1.7,5.46)-- (0.72,4.08);
\draw (0.72,4.38)-- (0.72,4.08);
\draw (0.72,4.08)-- (1.56,5.46);
\draw (1.56,5.78)-- (1.56,5.46);
\draw [rotate around={-17.53:(0.13,5.38)}] (0.13,5.38) ellipse (0.03cm and 0.03cm);
\draw (0.42,4.8)-- (0.74,4.8);
\draw (0.74,4.8)-- (0.58,5.04);
\draw (0.58,5.04)-- (0.96,5.1);
\draw (0.96,5.1)-- (0.78,5.32);
\draw (0.78,5.32)-- (1.14,5.42);
\draw (1.14,5.42)-- (0.92,5.56);
\draw (0.92,5.56)-- (1.26,5.64);
\end{tikzpicture})=
J_f(\begin{tikzpicture}[line cap=round,line join=round,>=triangle 45,x=0.11711167869707381cm,y=0.11711167869707381cm]
\clip(2.62,3.92) rectangle (6.04,6.09);
\draw (2.66,5.92)-- (5.92,5.92);
\draw (2.66,5.92)-- (3.48,4.52);
\draw (5.92,5.92)-- (3.48,4.52);
\draw (2.66,5.92)-- (2.66,5.44);
\draw (3.48,4.52)-- (3.5,4.02);
\draw (5.92,5.92)-- (5.94,5.46);
\draw (2.66,5.44)-- (3.5,4.02);
\draw (3.5,4.02)-- (5.94,5.46);
\draw (2.95,5.43)-- (5.09,5.44);
\end{tikzpicture})$.
\end{claim}

\begin{proof} \cite[\S10]{Bol77}
(a) 
Mark all the vertices of the polygons
$P,
P_1, \dots,
P_n$.
Marked points split the sides of the polygons into smaller
segments called \emph{edges}, see Figure~\ref{additivity}.
Each side is the sum of the edges, into which the side splits.
Thus the algebraic sum of the \emph{sides} can be replaced by the algebraic sum of the \emph{edges} in the definition of the additive invariant.

Consider the sum
$J_f(P_1)+\dots+J_f(P_n)$.
It equals to the algebraic sum of all the edges parallel to $f$ of all the polygons
$P_1,
P_2, \dots,
P_n$.

Take an edge inside the polygon $P$
(possibly excluding the endpoints).
Then there are two polygons $P_i$ and
$P_j$ bordering upon the edge from opposite sides.
The contributions of the edge to
$J_f(P_i)$ and
$J_f(P_j)$ have opposite signs and hence cancel in the sum
$J_f(P_1)+\dots+
J_f(P_n)$.

Take an edge in the boundary of the polygon
$P$. 
Then there is a unique polygon
$P_i$ bordering upon the edge.
The polygon $P_i$
borders upon the edge from the same side as the polygon
$P$.
Thus the contributions of the edge to $J_f(P_i)$
and $J_f(P)$ are the same.
Hence $J_f(P_1)+\dots+J_f(P_n)=
J_f(P)$.

\begin{figure}[htb]
{\fontsize{9pt}{9pt}\selectfont\begin{tikzpicture}[line cap=round,line join=round,>=triangle 45,x=0.66cm,y=0.66cm]
\clip(-2.36,-2.52) rectangle (6.64,5.22);
\fill[color=ffzzzz,fill=ffzzzz,fill opacity=0.2] (2.5,1) -- (0,1) -- (-1,2) -- (-2.28,2) -- (-1.2,1) -- (0,-1) -- (1,-2) -- (4,-2) -- (6,0) -- (6,1) -- cycle;
\fill[color=ffzzzz,fill=ffzzzz,fill opacity=0.2] (2.5,1.52) -- (2.5,3.52) -- (1,3.52) -- (0,1.52) -- cycle;
\fill[color=ffzzzz,fill=ffzzzz,fill opacity=0.2] (6.5,2.52) -- (5.5,2.52) -- (4.5,3.52) -- (3,3.52) -- (3,1.52) -- (6.5,1.52) -- cycle;
\draw [->,line width=1.2pt,color=zzqqcc] (5,4.5) -- (-1,4.5);
\draw [->,line width=1.2pt,color=zzqqcc] (-1,2) -- (-2.28,2);
\draw [->,color=wwttqq] (-2.28,2) -- (-1.2,1);
\draw [->,color=wwttqq] (-1.2,1) -- (0,-1);
\draw [->,color=wwttqq] (0,-1) -- (1,-2);
\draw [->,line width=1.2pt,color=zzqqcc] (1,-2) -- (4,-2);
\draw [->,color=wwttqq] (4,-2) -- (6,0);
\draw [->,color=wwttqq] (0,1) -- (-1,2);
\draw [->,color=wwttqq] (6,0) -- (6,1);
\draw [->,line width=1.2pt,color=zzqqcc] (6,1) -- (2.5,1);
\draw [->,line width=1.2pt,color=zzqqcc] (2.5,1) -- (0,1);
\draw [->,color=wwttqq] (2.5,1.52) -- (2.5,3.52);
\draw [->,line width=1.2pt,color=zzqqcc] (2.5,3.52) -- (1,3.52);
\draw [->,color=wwttqq] (1,3.52) -- (0,1.52);
\draw [->,line width=1.2pt,color=zzqqcc] (0,1.52) -- (2.5,1.52);
\draw [->,line width=1.2pt,color=zzqqcc] (3,1.52) -- (6.5,1.52);
\draw [->,color=wwttqq] (6.5,1.52) -- (6.5,2.52);
\draw [->,line width=1.2pt,color=zzqqcc] (6.5,2.52) -- (5.5,2.52);
\draw [->,color=wwttqq] (5.5,2.52) -- (4.5,3.52);
\draw [->,line width=1.2pt,color=zzqqcc] (4.5,3.52) -- (3,3.52);
\draw [->,color=wwttqq] (3,3.52) -- (3,1.52);
\draw[color=zzqqcc] (2,4.8) node {$f$};
\draw[color=zzqqcc] (-1.62,2.42) node {+};
\draw[color=zzqqcc] (2.62,-2.22) node {$-$};
\draw[color=ffzzzz] (2.44,-0.28) node {$P_1$};
\draw[color=ffzzzz] (1.7,2.66) node {$P_3$};
\draw[color=ffzzzz] (4.52,2.62) node {$P_2$};
\draw[color=zzqqcc] (4.36,0.68) node {+};
\draw[color=zzqqcc] (1.32,0.66) node {+};
\draw[color=zzqqcc] (1.8,3.92) node {+};
\draw[color=zzqqcc] (1.4,2.02) node {$-$};
\draw[color=zzqqcc] (4.8,2.02) node {$-$};
\draw[color=zzqqcc] (6.06,2.9) node {+};
\draw[color=zzqqcc] (3.8,3.94) node {+};
\end{tikzpicture}}
\caption{Additivity of the invariant:
$J_f(P)=J_f(P_1)+J_f(P_2)+J_f(P_3)$}
\label{additivity}
\end{figure}
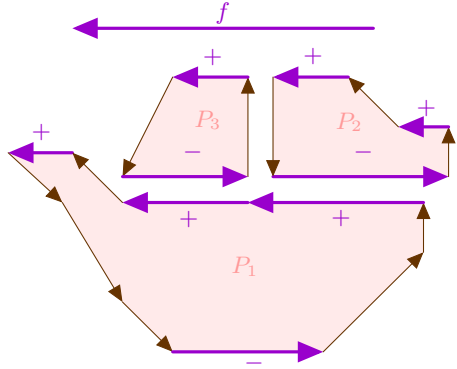

(b) Suppose that the cake 
$\begin{tikzpicture}[line cap=round,line join=round,>=triangle 45,x=0.11452372190810167cm,y=0.11452372190810167cm]
\clip(-1.75,3.88) rectangle (1.74,5.89);
\draw (-1.72,5.78)-- (0.72,4.38);
\draw (0.72,4.38)-- (1.56,5.78);
\draw (-1.72,5.78)-- (1.56,5.78);
\draw (-1.72,5.78)-- (-1.7,5.46);
\draw (-1.7,5.46)-- (0.72,4.08);
\draw (0.72,4.38)-- (0.72,4.08);
\draw (0.72,4.08)-- (1.56,5.46);
\draw (1.56,5.78)-- (1.56,5.46);
\draw [rotate around={-17.53:(0.13,5.38)}] (0.13,5.38) ellipse (0.03cm and 0.03cm);
\draw (0.42,4.8)-- (0.74,4.8);
\draw (0.74,4.8)-- (0.58,5.04);
\draw (0.58,5.04)-- (0.96,5.1);
\draw (0.96,5.1)-- (0.78,5.32);
\draw (0.78,5.32)-- (1.14,5.42);
\draw (1.14,5.42)-- (0.92,5.56);
\draw (0.92,5.56)-- (1.26,5.64);
\end{tikzpicture} $ is cut into polygonal pieces
$P_1,\dots,P_n$,
the box $\begin{tikzpicture}[line cap=round,line join=round,>=triangle 45,x=0.11711167869707381cm,y=0.11711167869707381cm]
\clip(2.62,3.92) rectangle (6.04,6.09);
\draw (2.66,5.92)-- (5.92,5.92);
\draw (2.66,5.92)-- (3.48,4.52);
\draw (5.92,5.92)-- (3.48,4.52);
\draw (2.66,5.92)-- (2.66,5.44);
\draw (3.48,4.52)-- (3.5,4.02);
\draw (5.92,5.92)-- (5.94,5.46);
\draw (2.66,5.44)-- (3.5,4.02);
\draw (3.5,4.02)-- (5.94,5.46);
\draw (2.95,5.43)-- (5.09,5.44);
\end{tikzpicture} $ is
cut into polygonal pieces
$P'_1,\dots,P'_n$,
and the piece $P'_i$ is a parallel translation of the piece $P_i$ for each $i=1,\dots,n$. Clearly,
$J_f(P_i)=J_f(P'_i)$.
Thus by assertion~(a) we have
$$
J_f(\begin{tikzpicture}[line cap=round,line join=round,>=triangle 45,x=0.11452372190810167cm,y=0.11452372190810167cm]
\clip(-1.75,3.88) rectangle (1.74,5.89);
\draw (-1.72,5.78)-- (0.72,4.38);
\draw (0.72,4.38)-- (1.56,5.78);
\draw (-1.72,5.78)-- (1.56,5.78);
\draw (-1.72,5.78)-- (-1.7,5.46);
\draw (-1.7,5.46)-- (0.72,4.08);
\draw (0.72,4.38)-- (0.72,4.08);
\draw (0.72,4.08)-- (1.56,5.46);
\draw (1.56,5.78)-- (1.56,5.46);
\draw [rotate around={-17.53:(0.13,5.38)}] (0.13,5.38) ellipse (0.03cm and 0.03cm);
\draw (0.42,4.8)-- (0.74,4.8);
\draw (0.74,4.8)-- (0.58,5.04);
\draw (0.58,5.04)-- (0.96,5.1);
\draw (0.96,5.1)-- (0.78,5.32);
\draw (0.78,5.32)-- (1.14,5.42);
\draw (1.14,5.42)-- (0.92,5.56);
\draw (0.92,5.56)-- (1.26,5.64);
\end{tikzpicture} )=
J_f(P_1)+\dots+J_f(P_n)=
J_f(P'_1)+\dots+J_f(P'_n)=
J_f(\begin{tikzpicture}[line cap=round,line join=round,>=triangle 45,x=0.11711167869707381cm,y=0.11711167869707381cm]
\clip(2.62,3.92) rectangle (6.04,6.09);
\draw (2.66,5.92)-- (5.92,5.92);
\draw (2.66,5.92)-- (3.48,4.52);
\draw (5.92,5.92)-- (3.48,4.52);
\draw (2.66,5.92)-- (2.66,5.44);
\draw (3.48,4.52)-- (3.5,4.02);
\draw (5.92,5.92)-- (5.94,5.46);
\draw (2.66,5.44)-- (3.5,4.02);
\draw (3.5,4.02)-- (5.94,5.46);
\draw (2.95,5.43)-- (5.09,5.44);
\end{tikzpicture} ).
$$
\end{proof}

To explain our idea let us prove an assertion close to Theorem~\ref{th1}. This assertion also shows that the Bolyai--Gerwien theorem is not necessarily true if only parallel translations of the pieces are permitted.

\begin{claim} \label{cl0} It is impossible to cut a cake into polygonal pieces and put the pieces into a box all together by means of only parallel translations.
\end{claim}

\begin{proof}
Assume that the cake is cut as required.
Then $J_{\mathbf{AB}}(\begin{tikzpicture}[line cap=round,line join=round,>=triangle 45,x=0.11452372190810167cm,y=0.11452372190810167cm]
\clip(-1.75,3.88) rectangle (1.74,5.89);
\draw (-1.72,5.78)-- (0.72,4.38);
\draw (0.72,4.38)-- (1.56,5.78);
\draw (-1.72,5.78)-- (1.56,5.78);
\draw (-1.72,5.78)-- (-1.7,5.46);
\draw (-1.7,5.46)-- (0.72,4.08);
\draw (0.72,4.38)-- (0.72,4.08);
\draw (0.72,4.08)-- (1.56,5.46);
\draw (1.56,5.78)-- (1.56,5.46);
\draw [rotate around={-17.53:(0.13,5.38)}] (0.13,5.38) ellipse (0.03cm and 0.03cm);
\draw (0.42,4.8)-- (0.74,4.8);
\draw (0.74,4.8)-- (0.58,5.04);
\draw (0.58,5.04)-- (0.96,5.1);
\draw (0.96,5.1)-- (0.78,5.32);
\draw (0.78,5.32)-- (1.14,5.42);
\draw (1.14,5.42)-- (0.92,5.56);
\draw (0.92,5.56)-- (1.26,5.64);
\end{tikzpicture} )=AB$. By Claim~\ref{cl1}(b)
$J_{\mathbf{AB}}(\begin{tikzpicture}[line cap=round,line join=round,>=triangle 45,x=0.11711167869707381cm,y=0.11711167869707381cm]
\clip(2.62,3.92) rectangle (6.04,6.09);
\draw (2.66,5.92)-- (5.92,5.92);
\draw (2.66,5.92)-- (3.48,4.52);
\draw (5.92,5.92)-- (3.48,4.52);
\draw (2.66,5.92)-- (2.66,5.44);
\draw (3.48,4.52)-- (3.5,4.02);
\draw (5.92,5.92)-- (5.94,5.46);
\draw (2.66,5.44)-- (3.5,4.02);
\draw (3.5,4.02)-- (5.94,5.46);
\draw (2.95,5.43)-- (5.09,5.44);
\end{tikzpicture} )=
J_{\mathbf{AB}}(\begin{tikzpicture}[line cap=round,line join=round,>=triangle 45,x=0.11452372190810167cm,y=0.11452372190810167cm]
\clip(-1.75,3.88) rectangle (1.74,5.89);
\draw (-1.72,5.78)-- (0.72,4.38);
\draw (0.72,4.38)-- (1.56,5.78);
\draw (-1.72,5.78)-- (1.56,5.78);
\draw (-1.72,5.78)-- (-1.7,5.46);
\draw (-1.7,5.46)-- (0.72,4.08);
\draw (0.72,4.38)-- (0.72,4.08);
\draw (0.72,4.08)-- (1.56,5.46);
\draw (1.56,5.78)-- (1.56,5.46);
\draw [rotate around={-17.53:(0.13,5.38)}] (0.13,5.38) ellipse (0.03cm and 0.03cm);
\draw (0.42,4.8)-- (0.74,4.8);
\draw (0.74,4.8)-- (0.58,5.04);
\draw (0.58,5.04)-- (0.96,5.1);
\draw (0.96,5.1)-- (0.78,5.32);
\draw (0.78,5.32)-- (1.14,5.42);
\draw (1.14,5.42)-- (0.92,5.56);
\draw (0.92,5.56)-- (1.26,5.64);
\end{tikzpicture} )$.
Thus $J_{\mathbf{AB}}(\begin{tikzpicture}[line cap=round,line join=round,>=triangle 45,x=0.11711167869707381cm,y=0.11711167869707381cm]
\clip(2.62,3.92) rectangle (6.04,6.09);
\draw (2.66,5.92)-- (5.92,5.92);
\draw (2.66,5.92)-- (3.48,4.52);
\draw (5.92,5.92)-- (3.48,4.52);
\draw (2.66,5.92)-- (2.66,5.44);
\draw (3.48,4.52)-- (3.5,4.02);
\draw (5.92,5.92)-- (5.94,5.46);
\draw (2.66,5.44)-- (3.5,4.02);
\draw (3.5,4.02)-- (5.94,5.46);
\draw (2.95,5.43)-- (5.09,5.44);
\end{tikzpicture} )\ne0$. Hence the box has a side parallel to $AB$. Analogously the box has $2$ sides parallel to $BC$ and $CA$. Thus the box is either a parallel translation of the cake or central-symmetric to the cake. Since the box is a mirror image of the cake it follows that the cake is mirror-symmetric, i.e., isosceles. This contradiction proves the claim.
\end{proof}

For the proof of Theorem~\ref{th1} we need to generalize the invariant $J_f$ to make it invariant under a rotation.
Denote by $\mathbf{XY}$ the directed line passing through points $X\ne Y$ and directed from $X$ to $Y$.

\begin{definition*}[of the generalized additive invariant]
Let $f({\mathbf{XY}})$ be a function on the set of all directed lines ${\mathbf{XY}}$ in the plane satisfying the property:
$f({\mathbf{XY}})=-f({\mathbf{YX}})$. Let  $P=A_1A_2\dots A_n$ be a polygon whose vertices are enumerated counterclockwise.
Set
$$J_f(P)=f(\mathbf{A}_1\mathbf{A}_2)A_1A_2+f(\mathbf{A}_2\mathbf{A}_3)A_2A_3+\dots+f(\mathbf{A}_n\mathbf{A}_1)A_nA_1.$$
\end{definition*}



We say that a rotation $R$ \emph{preserves} the function $f$ if for each directed line ${\mathbf{XY}}$
we have $f({R(\mathbf{XY}}))=f({\mathbf{XY}})$. If $R$ preserves $f$ then $R$ \emph{preserves} the invariant $J_f$, i.e., $J_f(R(P))=J_f(P)$ for any polygon $P$.
The following claim is proved analogously to Claim~\ref{cl1}.

\begin{claim}\label{cl2} Suppose that a cake 
\begin{tikzpicture}[line cap=round,line join=round,>=triangle 45,x=0.11452372190810167cm,y=0.11452372190810167cm]
\clip(-1.75,3.88) rectangle (1.74,5.89);
\draw (-1.72,5.78)-- (0.72,4.38);
\draw (0.72,4.38)-- (1.56,5.78);
\draw (-1.72,5.78)-- (1.56,5.78);
\draw (-1.72,5.78)-- (-1.7,5.46);
\draw (-1.7,5.46)-- (0.72,4.08);
\draw (0.72,4.38)-- (0.72,4.08);
\draw (0.72,4.08)-- (1.56,5.46);
\draw (1.56,5.78)-- (1.56,5.46);
\draw [rotate around={-17.53:(0.13,5.38)}] (0.13,5.38) ellipse (0.03cm and 0.03cm);
\draw (0.42,4.8)-- (0.74,4.8);
\draw (0.74,4.8)-- (0.58,5.04);
\draw (0.58,5.04)-- (0.96,5.1);
\draw (0.96,5.1)-- (0.78,5.32);
\draw (0.78,5.32)-- (1.14,5.42);
\draw (1.14,5.42)-- (0.92,5.56);
\draw (0.92,5.56)-- (1.26,5.64);
\end{tikzpicture}  is cut into several polygons and the pieces are put into a box 
\begin{tikzpicture}[line cap=round,line join=round,>=triangle 45,x=0.11711167869707381cm,y=0.11711167869707381cm]
\clip(2.62,3.92) rectangle (6.04,6.09);
\draw (2.66,5.92)-- (5.92,5.92);
\draw (2.66,5.92)-- (3.48,4.52);
\draw (5.92,5.92)-- (3.48,4.52);
\draw (2.66,5.92)-- (2.66,5.44);
\draw (3.48,4.52)-- (3.5,4.02);
\draw (5.92,5.92)-- (5.94,5.46);
\draw (2.66,5.44)-- (3.5,4.02);
\draw (3.5,4.02)-- (5.94,5.46);
\draw (2.95,5.43)-- (5.09,5.44);
\end{tikzpicture}  all together by means of rotations preserving the function $f$.
Then
$J_f(\begin{tikzpicture}[line cap=round,line join=round,>=triangle 45,x=0.11452372190810167cm,y=0.11452372190810167cm]
\clip(-1.75,3.88) rectangle (1.74,5.89);
\draw (-1.72,5.78)-- (0.72,4.38);
\draw (0.72,4.38)-- (1.56,5.78);
\draw (-1.72,5.78)-- (1.56,5.78);
\draw (-1.72,5.78)-- (-1.7,5.46);
\draw (-1.7,5.46)-- (0.72,4.08);
\draw (0.72,4.38)-- (0.72,4.08);
\draw (0.72,4.08)-- (1.56,5.46);
\draw (1.56,5.78)-- (1.56,5.46);
\draw [rotate around={-17.53:(0.13,5.38)}] (0.13,5.38) ellipse (0.03cm and 0.03cm);
\draw (0.42,4.8)-- (0.74,4.8);
\draw (0.74,4.8)-- (0.58,5.04);
\draw (0.58,5.04)-- (0.96,5.1);
\draw (0.96,5.1)-- (0.78,5.32);
\draw (0.78,5.32)-- (1.14,5.42);
\draw (1.14,5.42)-- (0.92,5.56);
\draw (0.92,5.56)-- (1.26,5.64);
\end{tikzpicture} )=
J_f(\begin{tikzpicture}[line cap=round,line join=round,>=triangle 45,x=0.11711167869707381cm,y=0.11711167869707381cm]
\clip(2.62,3.92) rectangle (6.04,6.09);
\draw (2.66,5.92)-- (5.92,5.92);
\draw (2.66,5.92)-- (3.48,4.52);
\draw (5.92,5.92)-- (3.48,4.52);
\draw (2.66,5.92)-- (2.66,5.44);
\draw (3.48,4.52)-- (3.5,4.02);
\draw (5.92,5.92)-- (5.94,5.46);
\draw (2.66,5.44)-- (3.5,4.02);
\draw (3.5,4.02)-- (5.94,5.46);
\draw (2.95,5.43)-- (5.09,5.44);
\end{tikzpicture} )$.
\end{claim}


\section{Proof of the theorem}

We are ready to prove Theorem~\ref{th1}.



\begin{figure}[htb]
{\fontsize{9pt}{9pt}\selectfont\begin{tikzpicture}[line cap=round,line join=round,>=triangle 45,x=1.0cm,y=1.0cm]
\clip(0.5,-3) rectangle (12.5,3);
\fill[color=ffzzzz,fill=ffzzzz,fill opacity=0.2] (0.74,0.3) -- (3.04,2.92) -- (5.31,0.28) -- (4.24,0.89) -- (3.02,1.1) -- (1.8,0.89) -- cycle;
\fill[line width=2pt,color=wwttqq,fill=wwttqq,fill opacity=0.4] (0.74,0.3) -- (1.8,0.89) -- (3.02,1.1) -- (4.24,0.89) -- (5.31,0.28) -- (6.12,-0.66) -- cycle;
\fill[color=ffzzzz,fill=ffzzzz,fill opacity=0.2] (12.4,0.3) -- (10.1,2.92) -- (7.82,0.28) -- (8.9,0.89) -- (10.11,1.1) -- (11.33,0.89) -- cycle;
\fill[line width=2pt,color=wwttqq,fill=wwttqq,fill opacity=0.4] (12.4,0.3) -- (11.33,0.89) -- (10.11,1.1) -- (8.9,0.89) -- (7.82,0.28) -- (7.02,-0.66) -- cycle;
\draw [shift={(3.01,-2.52)},color=wwttqq,fill=wwttqq,fill opacity=0.4] (0,0) -- (70.2:1.4) arc (70.2:89.86:1.4) -- cycle;
\draw [color=wwttqq] (0.74,0.3)-- (1.8,0.89);
\draw [color=wwttqq] (1.8,0.89)-- (3.02,1.1);
\draw [color=wwttqq] (3.02,1.1)-- (4.24,0.89);
\draw [color=wwttqq] (4.24,0.89)-- (5.31,0.28);
\draw [color=wwttqq] (5.31,0.28)-- (6.12,-0.66);
\draw [color=wwttqq] (6.12,-0.66)-- (0.74,0.3);
\draw [color=wwttqq] (10.1,2.92)-- (7.82,0.28);
\draw [color=wwttqq] (12.4,0.3)-- (11.33,0.89);
\draw [color=wwttqq] (11.33,0.89)-- (10.11,1.1);
\draw [color=wwttqq] (10.11,1.1)-- (8.9,0.89);
\draw [color=wwttqq] (8.9,0.89)-- (7.82,0.28);
\draw [color=wwttqq] (7.82,0.28)-- (7.02,-0.66);
\draw [color=wwttqq] (7.02,-0.66)-- (12.4,0.3);
\draw [shift={(3.01,-2.52)},->,color=wwttqq] (70.2:1.4) arc (70.2:89.86:1.4);
\draw [->,color=wwttqq] (5.31,0.28) -- (3.04,2.92);
\draw [->,color=wwttqq] (5.31,0.28) -- (4.24,0.89);
\draw [->,color=wwttqq] (4.24,0.89) -- (3.02,1.1);
\draw [->,color=wwttqq] (3.02,1.1) -- (1.8,0.89);
\draw [->,color=wwttqq] (1.8,0.89) -- (0.74,0.3);
\draw [->,color=wwttqq] (3.04,2.92) -- (0.74,0.3);
\draw [->,color=wwttqq] (7.82,0.28) -- (7.02,-0.66);
\draw [->,color=wwttqq] (8.9,0.89) -- (7.82,0.28);
\draw [->,color=wwttqq] (10.11,1.1) -- (8.9,0.89);
\draw [->,color=wwttqq] (11.33,0.89) -- (10.11,1.1);
\draw [->,color=wwttqq] (12.4,0.3) -- (11.33,0.89);
\draw [->,color=wwttqq] (12.4,0.3) -- (10.1,2.92);
\draw [color=wwttqq](3.14,-0.52) node[anchor=north west] {$R$};
\fill [color=wwttqq] (0.74,0.3) circle (1.0pt);
\draw[color=wwttqq] (0.8,0.05) node {$C$};
\fill [color=wwttqq] (6.12,-0.66) circle (1.0pt);
\draw[color=wwttqq] (6.25,-0.41) node {$A$};
\fill [color=wwttqq] (3.04,2.92) circle (1.0pt);
\draw[color=wwttqq] (3.42,2.88) node {$B$};
\fill [color=wwttqq] (10.1,2.92) circle (1.0pt);
\draw[color=wwttqq] (9.74,2.88) node {$B'$};
\fill [color=wwttqq] (12.4,0.3) circle (1.0pt);
\draw[color=wwttqq] (12.2,0.07) node {$C'$};
\fill [color=wwttqq] (7.02,-0.66) circle (1.0pt);
\draw[color=wwttqq] (6.83,-0.41) node {$A'$};
\fill [color=wwttqq] (3.01,-2.52) circle (1.0pt);
\draw[color=wwttqq] (3.2,-2.43) node {$O$};
\draw[color=wwttqq] (3.16,-1.51) node {$\phi$};
\draw[color=wwttqq] (4.75,0.15) node {$R(\mathbf{AB})$};
\draw[color=wwttqq] (3.57,0.64) node {$R^2(\mathbf{AB})$};
\fill [color=wwttqq] (10.12,-2.52) circle (1.0pt);
\draw[color=wwttqq] (10.42,-2.27) node {$O$};
\end{tikzpicture}}
\caption{The rotation $R$}
\label{construction}
\end{figure}
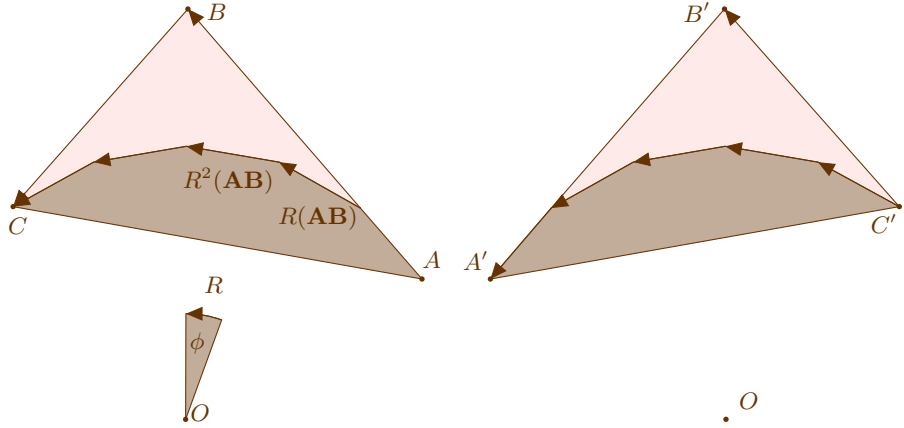

Let the cake \begin{tikzpicture}[line cap=round,line join=round,>=triangle 45,x=0.11452372190810167cm,y=0.11452372190810167cm]
\clip(-1.75,3.88) rectangle (1.74,5.89);
\draw (-1.72,5.78)-- (0.72,4.38);
\draw (0.72,4.38)-- (1.56,5.78);
\draw (-1.72,5.78)-- (1.56,5.78);
\draw (-1.72,5.78)-- (-1.7,5.46);
\draw (-1.7,5.46)-- (0.72,4.08);
\draw (0.72,4.38)-- (0.72,4.08);
\draw (0.72,4.08)-- (1.56,5.46);
\draw (1.56,5.78)-- (1.56,5.46);
\draw [rotate around={-17.53:(0.13,5.38)}] (0.13,5.38) ellipse (0.03cm and 0.03cm);
\draw (0.42,4.8)-- (0.74,4.8);
\draw (0.74,4.8)-- (0.58,5.04);
\draw (0.58,5.04)-- (0.96,5.1);
\draw (0.96,5.1)-- (0.78,5.32);
\draw (0.78,5.32)-- (1.14,5.42);
\draw (1.14,5.42)-- (0.92,5.56);
\draw (0.92,5.56)-- (1.26,5.64);
\end{tikzpicture}  be nicely cut into $2$ pieces. We may assume that while packing into the box \begin{tikzpicture}[line cap=round,line join=round,>=triangle 45,x=0.11711167869707381cm,y=0.11711167869707381cm]
\clip(2.62,3.92) rectangle (6.04,6.09);
\draw (2.66,5.92)-- (5.92,5.92);
\draw (2.66,5.92)-- (3.48,4.52);
\draw (5.92,5.92)-- (3.48,4.52);
\draw (2.66,5.92)-- (2.66,5.44);
\draw (3.48,4.52)-- (3.5,4.02);
\draw (5.92,5.92)-- (5.94,5.46);
\draw (2.66,5.44)-- (3.5,4.02);
\draw (3.5,4.02)-- (5.94,5.46);
\draw (2.95,5.43)-- (5.09,5.44);
\end{tikzpicture}  one piece remains fixed and the other piece is moved
by a rotation or parallel translation $R$.
By Claim~\ref{cl0} it follows that $R$ is not a parallel translation. Thus $R$ is a rotation through an angle $\phi$ around a point~$O$, see Figure~\ref{construction}.

We say that an angle $\psi$ is a \emph{multiple} of $\phi$, if $\psi=k\phi+l\cdot 180^\circ$ for some integers $k$ and $l$.
Notation: $\phi\mathbin{|}\psi$. An angle $\psi$ is \emph{rational} if $\psi=k\cdot 180^\circ/l$ for some integers $k$ and $l$.

Our aim is to show that one of the following conditions holds:
\begin{itemize}
\item each of the angles $2\alpha$, $2\beta$, $2\gamma$ is a multiple of $\phi$ (as in Example~\ref{ex1}(e)); or
\item one of the angles $2\alpha$, $2\beta$, $2\gamma$ is a multiple of $\phi$ and $\phi$ is rational (as in Example~\ref{ex1}(a)).
\end{itemize}





The \emph{oriented angle} $\angle(KL,MN)$ between (undirected) lines
$KL$ and $MN$ is an angle of a counterclockwise rotation taking the line $KL$ to $MN$.
An oriented angle is defined up to a multiple of $180^\circ$.

\begin{claim}\label{cl6} If either $\phi$ is irrational or $O\not\in \mathbf{AB}$ then
$\phi\mathbin{|}\angle(AB,A'B')$. 
\end{claim}

\begin{proof}
Let us construct an invariant $J_f$ preserved by the rotation $R$. Denote by $R^1,R^2,R^3,\dots$ the iterations of the rotation $R$,
by $R^{-1},R^{-2},R^{-3},\dots$ the iterations of the inverse rotation $R^{-1}$, and by $R^{0}$ the identity map. Define a function $f$ on the set of all directed lines by the formula
$$
f(\mathbf{XY})=
\begin{cases}
1,  & \quad\text{if $\mathbf{XY}=R^k(\mathbf{AB})$ for some integer $k$;}\\
-1, & \quad\text{if $\mathbf{XY}=R^k(\mathbf{BA})$ for some integer $k$;}\\
0, & \quad\text{otherwise.}
\end{cases}
$$
For instance, in Figure~\ref{construction} $f(\mathbf{AB})=1$, $f(\mathbf{A'B'})=-1$, $f(\mathbf{AC})=0$.

Since either $\phi$ is irrational 
or $O\not\in \mathbf{AB}$ it follows that $R^k(\mathbf{AB})\ne R^l(\mathbf{BA})$ for all integers $k,l$.
Thus the function $f$ is well-defined.
Clearly, the rotation $R$ preserves the function $f$. Notice that the condition $f(\mathbf{XY})\ne0$ implies $\phi\mathbin{|}\angle({AB},XY)$.


Consider the invariants $J_{f}(\begin{tikzpicture}[line cap=round,line join=round,>=triangle 45,x=0.11452372190810167cm,y=0.11452372190810167cm]
\clip(-1.75,3.88) rectangle (1.74,5.89);
\draw (-1.72,5.78)-- (0.72,4.38);
\draw (0.72,4.38)-- (1.56,5.78);
\draw (-1.72,5.78)-- (1.56,5.78);
\draw (-1.72,5.78)-- (-1.7,5.46);
\draw (-1.7,5.46)-- (0.72,4.08);
\draw (0.72,4.38)-- (0.72,4.08);
\draw (0.72,4.08)-- (1.56,5.46);
\draw (1.56,5.78)-- (1.56,5.46);
\draw [rotate around={-17.53:(0.13,5.38)}] (0.13,5.38) ellipse (0.03cm and 0.03cm);
\draw (0.42,4.8)-- (0.74,4.8);
\draw (0.74,4.8)-- (0.58,5.04);
\draw (0.58,5.04)-- (0.96,5.1);
\draw (0.96,5.1)-- (0.78,5.32);
\draw (0.78,5.32)-- (1.14,5.42);
\draw (1.14,5.42)-- (0.92,5.56);
\draw (0.92,5.56)-- (1.26,5.64);
\end{tikzpicture} )$
and $J_{f}(\begin{tikzpicture}[line cap=round,line join=round,>=triangle 45,x=0.11711167869707381cm,y=0.11711167869707381cm]
\clip(2.62,3.92) rectangle (6.04,6.09);
\draw (2.66,5.92)-- (5.92,5.92);
\draw (2.66,5.92)-- (3.48,4.52);
\draw (5.92,5.92)-- (3.48,4.52);
\draw (2.66,5.92)-- (2.66,5.44);
\draw (3.48,4.52)-- (3.5,4.02);
\draw (5.92,5.92)-- (5.94,5.46);
\draw (2.66,5.44)-- (3.5,4.02);
\draw (3.5,4.02)-- (5.94,5.46);
\draw (2.95,5.43)-- (5.09,5.44);
\end{tikzpicture} )$.
Both of them are linear combinations of the sides $AB$, $BC$, $CA$ with coefficients $0, \pm1$. For instance, in Figure \ref{construction}
$J_{f}(\begin{tikzpicture}[line cap=round,line join=round,>=triangle 45,x=0.11452372190810167cm,y=0.11452372190810167cm]
\clip(-1.75,3.88) rectangle (1.74,5.89);
\draw (-1.72,5.78)-- (0.72,4.38);
\draw (0.72,4.38)-- (1.56,5.78);
\draw (-1.72,5.78)-- (1.56,5.78);
\draw (-1.72,5.78)-- (-1.7,5.46);
\draw (-1.7,5.46)-- (0.72,4.08);
\draw (0.72,4.38)-- (0.72,4.08);
\draw (0.72,4.08)-- (1.56,5.46);
\draw (1.56,5.78)-- (1.56,5.46);
\draw [rotate around={-17.53:(0.13,5.38)}] (0.13,5.38) ellipse (0.03cm and 0.03cm);
\draw (0.42,4.8)-- (0.74,4.8);
\draw (0.74,4.8)-- (0.58,5.04);
\draw (0.58,5.04)-- (0.96,5.1);
\draw (0.96,5.1)-- (0.78,5.32);
\draw (0.78,5.32)-- (1.14,5.42);
\draw (1.14,5.42)-- (0.92,5.56);
\draw (0.92,5.56)-- (1.26,5.64);
\end{tikzpicture} )=
J_{f}(\begin{tikzpicture}[line cap=round,line join=round,>=triangle 45,x=0.11711167869707381cm,y=0.11711167869707381cm]
\clip(2.62,3.92) rectangle (6.04,6.09);
\draw (2.66,5.92)-- (5.92,5.92);
\draw (2.66,5.92)-- (3.48,4.52);
\draw (5.92,5.92)-- (3.48,4.52);
\draw (2.66,5.92)-- (2.66,5.44);
\draw (3.48,4.52)-- (3.5,4.02);
\draw (5.92,5.92)-- (5.94,5.46);
\draw (2.66,5.44)-- (3.5,4.02);
\draw (3.5,4.02)-- (5.94,5.46);
\draw (2.95,5.43)-- (5.09,5.44);
\end{tikzpicture} )=AB+BC$.

Let us prove that at least one of the sides $AB$, $BC$, $CA$ has nonzero contribution to both $J_{f}(\begin{tikzpicture}[line cap=round,line join=round,>=triangle 45,x=0.11452372190810167cm,y=0.11452372190810167cm]
\clip(-1.75,3.88) rectangle (1.74,5.89);
\draw (-1.72,5.78)-- (0.72,4.38);
\draw (0.72,4.38)-- (1.56,5.78);
\draw (-1.72,5.78)-- (1.56,5.78);
\draw (-1.72,5.78)-- (-1.7,5.46);
\draw (-1.7,5.46)-- (0.72,4.08);
\draw (0.72,4.38)-- (0.72,4.08);
\draw (0.72,4.08)-- (1.56,5.46);
\draw (1.56,5.78)-- (1.56,5.46);
\draw [rotate around={-17.53:(0.13,5.38)}] (0.13,5.38) ellipse (0.03cm and 0.03cm);
\draw (0.42,4.8)-- (0.74,4.8);
\draw (0.74,4.8)-- (0.58,5.04);
\draw (0.58,5.04)-- (0.96,5.1);
\draw (0.96,5.1)-- (0.78,5.32);
\draw (0.78,5.32)-- (1.14,5.42);
\draw (1.14,5.42)-- (0.92,5.56);
\draw (0.92,5.56)-- (1.26,5.64);
\end{tikzpicture} )$ and $J_{f}(\begin{tikzpicture}[line cap=round,line join=round,>=triangle 45,x=0.11711167869707381cm,y=0.11711167869707381cm]
\clip(2.62,3.92) rectangle (6.04,6.09);
\draw (2.66,5.92)-- (5.92,5.92);
\draw (2.66,5.92)-- (3.48,4.52);
\draw (5.92,5.92)-- (3.48,4.52);
\draw (2.66,5.92)-- (2.66,5.44);
\draw (3.48,4.52)-- (3.5,4.02);
\draw (5.92,5.92)-- (5.94,5.46);
\draw (2.66,5.44)-- (3.5,4.02);
\draw (3.5,4.02)-- (5.94,5.46);
\draw (2.95,5.43)-- (5.09,5.44);
\end{tikzpicture} )$.
Indeed, otherwise the difference
$J_{f}(\begin{tikzpicture}[line cap=round,line join=round,>=triangle 45,x=0.11452372190810167cm,y=0.11452372190810167cm]
\clip(-1.75,3.88) rectangle (1.74,5.89);
\draw (-1.72,5.78)-- (0.72,4.38);
\draw (0.72,4.38)-- (1.56,5.78);
\draw (-1.72,5.78)-- (1.56,5.78);
\draw (-1.72,5.78)-- (-1.7,5.46);
\draw (-1.7,5.46)-- (0.72,4.08);
\draw (0.72,4.38)-- (0.72,4.08);
\draw (0.72,4.08)-- (1.56,5.46);
\draw (1.56,5.78)-- (1.56,5.46);
\draw [rotate around={-17.53:(0.13,5.38)}] (0.13,5.38) ellipse (0.03cm and 0.03cm);
\draw (0.42,4.8)-- (0.74,4.8);
\draw (0.74,4.8)-- (0.58,5.04);
\draw (0.58,5.04)-- (0.96,5.1);
\draw (0.96,5.1)-- (0.78,5.32);
\draw (0.78,5.32)-- (1.14,5.42);
\draw (1.14,5.42)-- (0.92,5.56);
\draw (0.92,5.56)-- (1.26,5.64);
\end{tikzpicture} )-J_{f}(\begin{tikzpicture}[line cap=round,line join=round,>=triangle 45,x=0.11711167869707381cm,y=0.11711167869707381cm]
\clip(2.62,3.92) rectangle (6.04,6.09);
\draw (2.66,5.92)-- (5.92,5.92);
\draw (2.66,5.92)-- (3.48,4.52);
\draw (5.92,5.92)-- (3.48,4.52);
\draw (2.66,5.92)-- (2.66,5.44);
\draw (3.48,4.52)-- (3.5,4.02);
\draw (5.92,5.92)-- (5.94,5.46);
\draw (2.66,5.44)-- (3.5,4.02);
\draw (3.5,4.02)-- (5.94,5.46);
\draw (2.95,5.43)-- (5.09,5.44);
\end{tikzpicture} )$
is also a linear combination of these sides 
with coefficients $0, \pm1$.
Such a linear combination is nonzero by the triangle inequality and the assumption that
the cake is not isosceles. On the other hand, by Claim~\ref{cl2} we have
$J_{f}(\begin{tikzpicture}[line cap=round,line join=round,>=triangle 45,x=0.11452372190810167cm,y=0.11452372190810167cm]
\clip(-1.75,3.88) rectangle (1.74,5.89);
\draw (-1.72,5.78)-- (0.72,4.38);
\draw (0.72,4.38)-- (1.56,5.78);
\draw (-1.72,5.78)-- (1.56,5.78);
\draw (-1.72,5.78)-- (-1.7,5.46);
\draw (-1.7,5.46)-- (0.72,4.08);
\draw (0.72,4.38)-- (0.72,4.08);
\draw (0.72,4.08)-- (1.56,5.46);
\draw (1.56,5.78)-- (1.56,5.46);
\draw [rotate around={-17.53:(0.13,5.38)}] (0.13,5.38) ellipse (0.03cm and 0.03cm);
\draw (0.42,4.8)-- (0.74,4.8);
\draw (0.74,4.8)-- (0.58,5.04);
\draw (0.58,5.04)-- (0.96,5.1);
\draw (0.96,5.1)-- (0.78,5.32);
\draw (0.78,5.32)-- (1.14,5.42);
\draw (1.14,5.42)-- (0.92,5.56);
\draw (0.92,5.56)-- (1.26,5.64);
\end{tikzpicture} )-J_{f}(\begin{tikzpicture}[line cap=round,line join=round,>=triangle 45,x=0.11711167869707381cm,y=0.11711167869707381cm]
\clip(2.62,3.92) rectangle (6.04,6.09);
\draw (2.66,5.92)-- (5.92,5.92);
\draw (2.66,5.92)-- (3.48,4.52);
\draw (5.92,5.92)-- (3.48,4.52);
\draw (2.66,5.92)-- (2.66,5.44);
\draw (3.48,4.52)-- (3.5,4.02);
\draw (5.92,5.92)-- (5.94,5.46);
\draw (2.66,5.44)-- (3.5,4.02);
\draw (3.5,4.02)-- (5.94,5.46);
\draw (2.95,5.43)-- (5.09,5.44);
\end{tikzpicture} )=0$, a contradiction.

So a side of the cake has nonzero contribution to both $J_{f}(\begin{tikzpicture}[line cap=round,line join=round,>=triangle 45,x=0.11452372190810167cm,y=0.11452372190810167cm]
\clip(-1.75,3.88) rectangle (1.74,5.89);
\draw (-1.72,5.78)-- (0.72,4.38);
\draw (0.72,4.38)-- (1.56,5.78);
\draw (-1.72,5.78)-- (1.56,5.78);
\draw (-1.72,5.78)-- (-1.7,5.46);
\draw (-1.7,5.46)-- (0.72,4.08);
\draw (0.72,4.38)-- (0.72,4.08);
\draw (0.72,4.08)-- (1.56,5.46);
\draw (1.56,5.78)-- (1.56,5.46);
\draw [rotate around={-17.53:(0.13,5.38)}] (0.13,5.38) ellipse (0.03cm and 0.03cm);
\draw (0.42,4.8)-- (0.74,4.8);
\draw (0.74,4.8)-- (0.58,5.04);
\draw (0.58,5.04)-- (0.96,5.1);
\draw (0.96,5.1)-- (0.78,5.32);
\draw (0.78,5.32)-- (1.14,5.42);
\draw (1.14,5.42)-- (0.92,5.56);
\draw (0.92,5.56)-- (1.26,5.64);
\end{tikzpicture} )$ and
$J_{f}(\begin{tikzpicture}[line cap=round,line join=round,>=triangle 45,x=0.11711167869707381cm,y=0.11711167869707381cm]
\clip(2.62,3.92) rectangle (6.04,6.09);
\draw (2.66,5.92)-- (5.92,5.92);
\draw (2.66,5.92)-- (3.48,4.52);
\draw (5.92,5.92)-- (3.48,4.52);
\draw (2.66,5.92)-- (2.66,5.44);
\draw (3.48,4.52)-- (3.5,4.02);
\draw (5.92,5.92)-- (5.94,5.46);
\draw (2.66,5.44)-- (3.5,4.02);
\draw (3.5,4.02)-- (5.94,5.46);
\draw (2.95,5.43)-- (5.09,5.44);
\end{tikzpicture} )$.

In case when the side is $AB$ we get $f(\mathbf{A'B'})\ne0$. Thus $\phi\mathbin{|}\angle({AB},{A}'{B}')$.

In case when the side is $BC$ we get $f(\mathbf{BC})\ne 0$ and $f(\mathbf{B'C'})\ne0$.
Thus $\phi\mathbin{|}\angle({AB},{BC})$ and $\phi\mathbin{|}\angle({AB},{B'C'})$. Then
$$
\angle ({AB},{A'B'})=\angle({AB},{B'C'})+\angle({B'C'},{A'B'})=
\angle({AB},{B'C'})+\angle({AB},{BC})
$$
is a multiple of $\phi$.

Case when the side is $AC$ is analogous.
\end{proof}






\begin{claim}\label{cl7} If either $\phi$ is irrational or $O\not\in \mathbf{AB},\mathbf{AC}$ then $\phi\mathbin{|}2\alpha$. 
\end{claim}

\begin{proof} By Claim~\ref{cl6} we have  $\phi\mathbin{|}\angle({AB},{A}'{B}')$ and $\phi\mathbin{|}\angle({AC},{A}'{C}')$. Thus
$$
2\alpha=\angle({AB},{AC})+\angle({A'C'},{A'B'})=\angle({AB},{A'B'})-\angle({AC},{A'C'})
$$
is a multiple of $\phi$.
\end{proof}

\begin{claim}\label{cl8} At least two lines containing the sides of the cake do not pass through the point $O$.
\end{claim}

\begin{proof} Assume the converse. Then $O$ is a vertex of the cake. Assume without loss of generality
that $O=A$ and $AB>AC$.

Let us show that $A'=A$. Indeed, consider all the angles with vertex at $A$ belonging to both pieces.
The sum of the angles is $\alpha$ and is preserved by the rotation $R$.
Thus the box has angle of size $\alpha$ at the point $A$. That is, $A'=A$.

Let us show that either $B'=B$ or $B'=R(B)$. Indeed, the point $B$ of the cake is the most distant from $A$. Thus no other points of the cake border upon $B$ after the rearrangement of the pieces. Thus the box has angle of size $\beta$ at the point $B$ or $R(B)$, depending on whether the piece containing $B$ is fixed or mobile. That is, either $B'=B$ or $B'=R(B)$.

In case $A'=A$ and $B'=B$ we get
$\begin{tikzpicture}[line cap=round,line join=round,>=triangle 45,x=0.11452372190810167cm,y=0.11452372190810167cm]
\clip(-1.75,3.88) rectangle (1.74,5.89);
\draw (-1.72,5.78)-- (0.72,4.38);
\draw (0.72,4.38)-- (1.56,5.78);
\draw (-1.72,5.78)-- (1.56,5.78);
\draw (-1.72,5.78)-- (-1.7,5.46);
\draw (-1.7,5.46)-- (0.72,4.08);
\draw (0.72,4.38)-- (0.72,4.08);
\draw (0.72,4.08)-- (1.56,5.46);
\draw (1.56,5.78)-- (1.56,5.46);
\draw [rotate around={-17.53:(0.13,5.38)}] (0.13,5.38) ellipse (0.03cm and 0.03cm);
\draw (0.42,4.8)-- (0.74,4.8);
\draw (0.74,4.8)-- (0.58,5.04);
\draw (0.58,5.04)-- (0.96,5.1);
\draw (0.96,5.1)-- (0.78,5.32);
\draw (0.78,5.32)-- (1.14,5.42);
\draw (1.14,5.42)-- (0.92,5.56);
\draw (0.92,5.56)-- (1.26,5.64);
\end{tikzpicture} \cap\begin{tikzpicture}[line cap=round,line join=round,>=triangle 45,x=0.11711167869707381cm,y=0.11711167869707381cm]
\clip(2.62,3.92) rectangle (6.04,6.09);
\draw (2.66,5.92)-- (5.92,5.92);
\draw (2.66,5.92)-- (3.48,4.52);
\draw (5.92,5.92)-- (3.48,4.52);
\draw (2.66,5.92)-- (2.66,5.44);
\draw (3.48,4.52)-- (3.5,4.02);
\draw (5.92,5.92)-- (5.94,5.46);
\draw (2.66,5.44)-- (3.5,4.02);
\draw (3.5,4.02)-- (5.94,5.46);
\draw (2.95,5.43)-- (5.09,5.44);
\end{tikzpicture} =A'B'$. Thus the fixed piece is empty.

In case $A'=A$ and $B'=R(B)$ we get
$R(\begin{tikzpicture}[line cap=round,line join=round,>=triangle 45,x=0.11452372190810167cm,y=0.11452372190810167cm]
\clip(-1.75,3.88) rectangle (1.74,5.89);
\draw (-1.72,5.78)-- (0.72,4.38);
\draw (0.72,4.38)-- (1.56,5.78);
\draw (-1.72,5.78)-- (1.56,5.78);
\draw (-1.72,5.78)-- (-1.7,5.46);
\draw (-1.7,5.46)-- (0.72,4.08);
\draw (0.72,4.38)-- (0.72,4.08);
\draw (0.72,4.08)-- (1.56,5.46);
\draw (1.56,5.78)-- (1.56,5.46);
\draw [rotate around={-17.53:(0.13,5.38)}] (0.13,5.38) ellipse (0.03cm and 0.03cm);
\draw (0.42,4.8)-- (0.74,4.8);
\draw (0.74,4.8)-- (0.58,5.04);
\draw (0.58,5.04)-- (0.96,5.1);
\draw (0.96,5.1)-- (0.78,5.32);
\draw (0.78,5.32)-- (1.14,5.42);
\draw (1.14,5.42)-- (0.92,5.56);
\draw (0.92,5.56)-- (1.26,5.64);
\end{tikzpicture} )\cap\begin{tikzpicture}[line cap=round,line join=round,>=triangle 45,x=0.11711167869707381cm,y=0.11711167869707381cm]
\clip(2.62,3.92) rectangle (6.04,6.09);
\draw (2.66,5.92)-- (5.92,5.92);
\draw (2.66,5.92)-- (3.48,4.52);
\draw (5.92,5.92)-- (3.48,4.52);
\draw (2.66,5.92)-- (2.66,5.44);
\draw (3.48,4.52)-- (3.5,4.02);
\draw (5.92,5.92)-- (5.94,5.46);
\draw (2.66,5.44)-- (3.5,4.02);
\draw (3.5,4.02)-- (5.94,5.46);
\draw (2.95,5.43)-- (5.09,5.44);
\end{tikzpicture} =A'B'$. Thus the mobile piece is empty.

This contradiction proves the claim.
\end{proof}

Now the proof of the theorem is concluded by consideration of the following $2$ cases.

Case 1: $\phi$ is irrational. By Claim~\ref{cl7} it follows that $\phi\mathbin{|}2\alpha, 2\beta, 2\gamma$. Then
$k\alpha+l\beta+m\gamma=0$ for some integers $k$, $l$ and $m$, not vanishing simultaneously.

Case 2: $\phi$ is rational. By Claim~\ref{cl8} we may assume without loss of generality that $O\not\in AB, AC$.
Then by Claim~\ref{cl7} $\phi\mathbin{|}2\alpha$. Thus $\alpha=k\cdot 180^\circ/l$ for some integers $k$ and $l$. So
$0=k\cdot 180^\circ-l\alpha =(k-l)\alpha+k\beta+k\gamma$.

\section{An open problem}


Let us state an open problem, cf. \cite{Zak}. 
We say that two similar nonisosceles triangles in the plane are {\it oriented oppositely} if one of them includes angles  $\alpha$, $\beta$, $\gamma$ clockwise, and another one counterclockwise, see Figure~\ref{fig1}.
For instance, the height from the right angle cuts a right triangle into two triangles similar to it but oriented oppositely.



\begin{problem} 
Does there exist a nonright and nonisosceles triangle which can be cut into triangles similar to it but oriented oppositely?
\end{problem}

Let us announce a partial result due to M.V.~Prasolov and the author. The proof is based on a similar invariant as above.


\begin{theorem}\label{th2} Suppose that a triangle with angles $\alpha$, $\beta$ and $\gamma$ can be dissected into triangles similar to it but oriented oppositely. Then $k\alpha+l\beta+m\gamma=0$ for some integers $k$, $l$ and $m$, not vanishing simultaneously.
\end{theorem}

\pagebreak

\section{Answers}\label{final}

\begin{proof}[Proof of Example~\ref{ex1}\textup{(a)--(d)}] A useful observation is: if the pieces are mirror-symmetric then the cutting is automatically nice. This allows not to take care about rearrangement of the pieces.

(a) Cut the cake along the median from the vertex of the right angle, see Figure~\ref{figsol}(a).

(b) Cut the cake along the line separating an angle of size $\beta$ from the angle $\alpha$, see Figure~\ref{figsol}(b).

(c) Cut the cake along the line separating an angle of size $\beta$ from the angle $\gamma$, see Figure~\ref{figsol}(c).

(d) Cut the cake along the line symmetric to
the side opposite to the angle $\gamma$ with respect to the bisector of the angle $\gamma$, see Figure~\ref{figsol}(d).
\end{proof}

\begin{figure}[htbp]
\begin{tabular}{cc}
{\fontsize{9pt}{9pt}\selectfont\begin{tikzpicture}[line cap=round,line join=round,>=triangle 45,x=1.0cm,y=1.0cm]
\clip(-0.54,-0.1) rectangle (4.5,2.14);
\draw [shift={(4.11,0.05)},color=wwttqq,fill=wwttqq,fill opacity=0.2] (0,0) -- (155.04:1.2) arc (155.04:180:1.2) -- cycle;
\draw [shift={(0.05,1.94)},color=wwttqq,fill=wwttqq,fill opacity=0.2] (0,0) -- (-90:0.7) arc (-90:-24.96:0.7) -- cycle;
\fill[color=wwttqq,fill=wwttqq,fill opacity=0.4] (0.05,1.94) -- (2.08,1) -- (0.05,0.05) -- cycle;
\fill[color=ffzzzz,fill=ffzzzz,fill opacity=0.2] (2.08,1) -- (4.11,0.05) -- (0.05,0.05) -- cycle;
\draw [shift={(0.05,0.05)},color=wwttqq,fill=wwttqq,fill opacity=0.2] (0,0) -- (24.96:0.7) arc (24.96:90:0.7) -- cycle;
\draw [shift={(0.05,0.05)},color=wwttqq,fill=wwttqq,fill opacity=0.2] (0,0) -- (0:1.2) arc (0:24.96:1.2) -- cycle;
\draw [color=wwttqq] (0.05,1.94)-- (2.08,1);
\draw [color=wwttqq] (1.1,1.55) -- (1.03,1.39);
\draw [color=wwttqq] (2.08,1)-- (0.05,0.05);
\draw [color=wwttqq] (0.05,0.05)-- (0.05,1.94);
\draw [color=wwttqq] (2.08,1)-- (4.11,0.05);
\draw [color=wwttqq] (3.13,0.6) -- (3.06,0.44);
\draw [color=wwttqq] (4.11,0.05)-- (0.05,0.05);
\draw [color=wwttqq] (0.05,0.05)-- (2.08,1);
\draw [dotted,color=wwttqq] (4.11,0.05)-- (2.08,1);
\draw [dotted,color=wwttqq] (2.08,1)-- (4.11,1.94);
\draw [dotted,color=wwttqq] (4.11,1.94)-- (4.11,0.05);
\draw [shift={(2.08,1)},color=wwttqq]  plot[domain=0.99:2.1,variable=\t]({1*0.77*cos(\t r)+0*0.77*sin(\t r)},{0*0.77*cos(\t r)+1*0.77*sin(\t r)});
\draw [->,color=wwttqq] (2.5,1.64) -- (2.66,1.56);
\draw[color=wwttqq] (3.16,0.22) node {$\beta$};
\draw[color=wwttqq] (0.28,1.52) node {$\gamma$};
\draw[color=wwttqq] (0.22,0.32) node {$\gamma$};
\draw[color=wwttqq] (0.9,0.24) node {$\beta$};
\end{tikzpicture}}
&
{\fontsize{9pt}{9pt}\selectfont\begin{tikzpicture}[line cap=round,line join=round,>=triangle 45,x=1.0cm,y=1.0cm]
\clip(0,-0.1) rectangle (5.02,1.98);
\draw [shift={(4.88,0.05)},color=wwttqq,fill=wwttqq,fill opacity=0.2] (0,0) -- (157.44:1.2) arc (157.44:180:1.2) -- cycle;
\fill[color=wwttqq,fill=wwttqq,fill opacity=0.4] (0.05,0.05) -- (0.74,1.77) -- (2.47,1.05) -- cycle;
\fill[color=ffzzzz,fill=ffzzzz,fill opacity=0.2] (0.05,0.05) -- (2.47,1.05) -- (4.88,0.05) -- cycle;
\draw [shift={(0.05,0.05)},color=wwttqq,fill=wwttqq,fill opacity=0.2] (0,0) -- (22.56:1.2) arc (22.56:68.28:1.2) -- cycle;
\draw [shift={(0.05,0.05)},color=wwttqq,fill=wwttqq,fill opacity=0.2] (0,0) -- (0:1.2) arc (0:22.56:1.2) -- cycle;
\draw [shift={(2.47,1.05)},color=wwttqq,fill=wwttqq,fill opacity=0.2] (0,0) -- (157.44:1.2) arc (157.44:202.56:1.2) -- cycle;
\draw [color=wwttqq] (0.05,0.05)-- (0.74,1.77);
\draw [color=wwttqq] (0.74,1.77)-- (2.47,1.05);
\draw [color=wwttqq] (2.47,1.05)-- (0.05,0.05);
\draw [color=wwttqq] (0.05,0.05)-- (2.47,1.05);
\draw [color=wwttqq] (2.47,1.05)-- (4.88,0.05);
\draw [color=wwttqq] (4.88,0.05)-- (0.05,0.05);
\draw [dotted,color=wwttqq] (4.88,0.05)-- (4.19,1.77);
\draw [dotted,color=wwttqq] (4.19,1.77)-- (2.47,1.05);
\draw [dotted,color=wwttqq] (2.47,1.05)-- (4.88,0.05);
\draw [shift={(0.05,0.05)},color=wwttqq] (22.56:1.2) arc (22.56:68.28:1.2);
\draw [shift={(0.05,0.05)},color=wwttqq] (22.56:1.1) arc (22.56:68.28:1.1);
\draw [shift={(2.47,1.05)},color=wwttqq] (157.44:1.2) arc (157.44:202.56:1.2);
\draw [shift={(2.47,1.05)},color=wwttqq] (157.44:1.1) arc (157.44:202.56:1.1);
\draw [shift={(2.47,1.05)},color=wwttqq]  plot[domain=0.89:2.35,variable=\t]({1*0.73*cos(\t r)+0*0.73*sin(\t r)},{0*0.73*cos(\t r)+1*0.73*sin(\t r)});
\draw [->,color=wwttqq] (2.92,1.62) -- (3.06,1.52);
\draw[color=wwttqq] (4.02,0.18) node {$\beta$};
\draw[color=wwttqq] (0.66,0.54) node {$2\beta$};
\draw[color=wwttqq] (1.04,0.18) node {$\beta$};
\draw[color=wwttqq] (1.82,1) node {$2\beta$};
\end{tikzpicture}}
\\[3pt]
a & b
\\[3pt]
{\fontsize{9pt}{9pt}\selectfont\begin{tikzpicture}[line cap=round,line join=round,>=triangle 45,x=1.0cm,y=1.0cm]
\clip(-0.26,-2.32) rectangle (5.22,1.98);
\draw [shift={(0,0)},color=wwttqq,fill=wwttqq,fill opacity=0.2] (0,0) -- (0:0.8) arc (0:53.62:0.8) -- cycle;
\draw [shift={(5.04,0)},color=wwttqq,fill=wwttqq,fill opacity=0.2] (0,0) -- (152.44:1.1) arc (152.44:180:1.1) -- cycle;
\fill[color=wwttqq,fill=wwttqq,fill opacity=0.4] (0,0) -- (1.4,1.9) -- (2.72,0) -- cycle;
\fill[color=ffzzzz,fill=ffzzzz,fill opacity=0.2] (1.4,1.9) -- (5.04,0) -- (2.72,0) -- cycle;
\draw [shift={(2.72,0)},color=wwttqq,fill=wwttqq,fill opacity=0.2] (0,0) -- (124.87:0.8) arc (124.87:180:0.8) -- cycle;
\draw [shift={(1.4,1.9)},color=wwttqq,fill=wwttqq,fill opacity=0.2] (0,0) -- (-55.13:1.1) arc (-55.13:-27.56:1.1) -- cycle;
\draw [shift={(0,0)},color=wwttqq] (0:0.8) arc (0:53.62:0.8);
\draw [shift={(0,0)},color=wwttqq] (0:0.7) arc (0:53.62:0.7);
\draw [color=wwttqq] (0,0)-- (1.4,1.9);
\draw [color=wwttqq] (1.4,1.9)-- (2.72,0);
\draw [color=wwttqq] (2.72,0)-- (0,0);
\draw [color=wwttqq] (1.4,1.9)-- (5.04,0);
\draw [color=wwttqq] (5.04,0)-- (2.72,0);
\draw [color=wwttqq] (2.72,0)-- (1.4,1.9);
\draw [dotted,color=wwttqq] (4.28,-2.23)-- (5.04,0);
\draw [dotted,color=wwttqq] (2.72,0)-- (4.28,-2.23);
\draw [shift={(1.34,-2.65)},color=wwttqq]  plot[domain=0.67:1.56,variable=\t]({1*2.22*cos(\t r)+0*2.22*sin(\t r)},{0*2.22*cos(\t r)+1*2.22*sin(\t r)});
\draw [->,color=wwttqq] (3.08,-1.28) -- (3.2,-1.46);
\draw [shift={(2.72,0)},color=wwttqq] (124.87:0.8) arc (124.87:180:0.8);
\draw [shift={(2.72,0)},color=wwttqq] (124.87:0.7) arc (124.87:180:0.7);
\draw[color=wwttqq] (0.44,0.12) node {$2\beta$};
\draw[color=wwttqq] (4.18,0.14) node {$\beta$};
\draw[color=wwttqq] (2.34,0.12) node {$2\beta$};
\draw[color=wwttqq] (2.12,1.28) node {$\beta$};
\end{tikzpicture}}
&
{\fontsize{9pt}{9pt}\selectfont\begin{tikzpicture}[line cap=round,line join=round,>=triangle 45,x=1.0cm,y=1.0cm]
\clip(-0.06,-1.12) rectangle (3.06,3.62);
\draw [shift={(0.4,0)},color=wwttqq,fill=wwttqq,fill opacity=0.2] (0,0) -- (0:0.6) arc (0:96.43:0.6) -- cycle;
\draw [shift={(3,0)},color=wwttqq,fill=wwttqq,fill opacity=0.2] (0,0) -- (130.18:0.7) arc (130.18:180:0.7) -- cycle;
\fill[color=wwttqq,fill=wwttqq,fill opacity=0.4] (3,0) -- (2.31,0.82) -- (1.53,0) -- cycle;
\fill[color=ffzzzz,fill=ffzzzz,fill opacity=0.2] (0,3.55) -- (0.4,0) -- (1.53,0) -- (2.31,0.82) -- cycle;
\draw [shift={(1.53,0)},color=wwttqq,fill=wwttqq,fill opacity=0.2] (0,0) -- (0:0.7) arc (0:46.61:0.7) -- cycle;
\draw [shift={(2.31,0.82)},color=wwttqq,fill=wwttqq,fill opacity=0.2] (0,0) -- (130.18:0.6) arc (130.18:226.61:0.6) -- cycle;
\draw [shift={(0.4,0)},color=wwttqq] (0:0.6) arc (0:96.43:0.6);
\draw [shift={(0.4,0)},color=wwttqq] (0:0.5) arc (0:96.43:0.5);
\draw [color=wwttqq] (3,0)-- (2.31,0.82);
\draw [color=wwttqq] (2.31,0.82)-- (1.53,0);
\draw [color=wwttqq] (1.53,0)-- (3,0);
\draw [color=wwttqq] (0,3.55)-- (0.4,0);
\draw [color=wwttqq] (0.4,0)-- (1.53,0);
\draw [color=wwttqq] (1.53,0)-- (2.31,0.82);
\draw [color=wwttqq] (2.31,0.82)-- (0,3.55);
\draw [dotted,color=wwttqq] (0.52,-1.07)-- (0.4,0);
\draw [dotted,color=wwttqq] (1.53,0)-- (0.52,-1.07);
\draw [shift={(1.99,-1.07)},color=wwttqq]  plot[domain=1.41:2.48,variable=\t]({1*0.68*cos(\t r)+0*0.68*sin(\t r)},{0*0.68*cos(\t r)+1*0.68*sin(\t r)});
\draw [->,color=wwttqq] (1.46,-0.65) -- (1.32,-0.84);
\draw [dash pattern=on 1pt off 3pt on 5pt off 4pt,color=zzqqcc] (0,3.55)-- (1.53,0);
\draw [shift={(2.31,0.82)},color=wwttqq] (130.18:0.6) arc (130.18:226.61:0.6);
\draw [shift={(2.31,0.82)},color=wwttqq] (130.18:0.5) arc (130.18:226.61:0.5);
\draw[color=wwttqq] (0.62,0.14) node {$2\beta$};
\draw[color=wwttqq] (2.6,0.14) node {$\beta$};
\draw[color=wwttqq] (2,0.14) node {$\beta$};
\draw[color=wwttqq] (2.04,0.78) node {$2\beta$};
\end{tikzpicture}}
\\[3pt]
c & d
\end{tabular}
\caption{Nice cuttings of cakes from Figure~4} 
\label{figsol}
\end{figure}
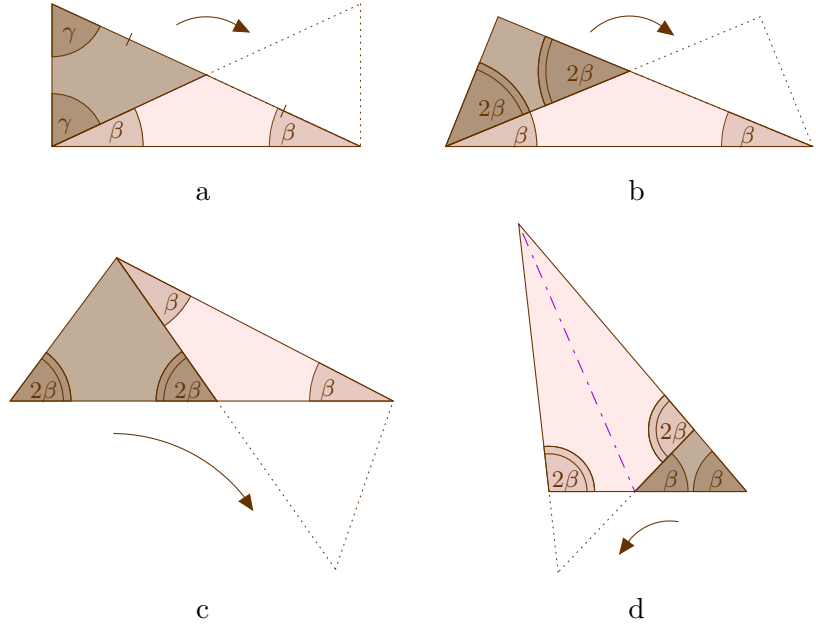

\begin{figure}[htbp]
{\fontsize{9pt}{9pt}\selectfont\begin{tikzpicture}[line cap=round,line join=round,>=triangle 45,x=1.0cm,y=1.0cm]
\clip(1.92,-2.47) rectangle (11.28,3.51);
\fill[color=ffzzzz,fill=ffzzzz,fill opacity=0.2] (11.04,1) -- (5.71,1) -- (5.19,1.62) -- (5.33,2.43) -- (6.04,2.84) -- cycle;
\fill[color=wwttqq,fill=wwttqq,fill opacity=0.45] (2.86,1) -- (5.71,1) -- (5.19,1.62) -- (5.33,2.43) -- cycle;
\draw [shift={(6.04,2.84)},color=ffzzzz,fill=ffzzzz,fill opacity=0.2] (0,0) -- (-149.95:0.3) arc (-149.95:-20.2:0.3) -- cycle;
\draw [shift={(5.33,2.43)},color=ffzzzz,fill=ffzzzz,fill opacity=0.2] (0,0) -- (-99.69:0.3) arc (-99.69:30.05:0.3) -- cycle;
\draw [shift={(5.19,1.62)},color=ffzzzz,fill=ffzzzz,fill opacity=0.2] (0,0) -- (-49.84:0.3) arc (-49.84:80.31:0.3) -- cycle;
\draw [shift={(5.71,1)},color=ffzzzz,fill=ffzzzz,fill opacity=0.2] (0,0) -- (0:0.3) arc (0:130.16:0.3) -- cycle;
\fill[color=ffzzzz,fill=ffzzzz,fill opacity=0.2] (11.04,0) -- (5.71,0) -- (5.19,-0.62) -- (5.33,-1.43) -- (6.04,-1.84) -- cycle;
\fill[color=wwttqq,fill=wwttqq,fill opacity=0.4] (2.86,0) -- (5.71,0) -- (5.19,-0.62) -- (5.33,-1.43) -- cycle;
\draw [shift={(5.71,0)},color=ffzzzz,fill=ffzzzz,fill opacity=0.2] (0,0) -- (-130.16:0.3) arc (-130.16:0:0.3) -- cycle;
\draw [shift={(5.19,-0.62)},color=ffzzzz,fill=ffzzzz,fill opacity=0.2] (0,0) -- (-80.31:0.3) arc (-80.31:49.84:0.3) -- cycle;
\draw [shift={(5.33,-1.43)},color=ffzzzz,fill=ffzzzz,fill opacity=0.2] (0,0) -- (-30.05:0.3) arc (-30.05:99.69:0.3) -- cycle;
\draw [shift={(6.04,-1.84)},color=ffzzzz,fill=ffzzzz,fill opacity=0.2] (0,0) -- (20.2:0.3) arc (20.2:149.95:0.3) -- cycle;
\draw [color=wwttqq] (11.04,1)-- (5.71,1);
\draw [color=wwttqq] (5.71,1)-- (5.19,1.62);
\draw [color=wwttqq] (5.38,1.25) -- (5.52,1.37);
\draw [color=wwttqq] (5.19,1.62)-- (5.33,2.43);
\draw [color=wwttqq] (5.17,2.04) -- (5.35,2.01);
\draw [color=wwttqq] (5.33,2.43)-- (6.04,2.84);
\draw [color=wwttqq] (5.64,2.71) -- (5.73,2.56);
\draw [color=wwttqq] (6.04,2.84)-- (11.04,1);
\draw [color=wwttqq] (2.86,1)-- (5.71,1);
\draw [color=wwttqq] (4.25,1.09) -- (4.25,0.91);
\draw [color=wwttqq] (4.32,1.09) -- (4.32,0.91);
\draw [color=wwttqq] (5.33,2.43)-- (2.86,1);
\draw [color=wwttqq] (4.17,1.65) -- (4.08,1.81);
\draw [color=wwttqq] (4.11,1.62) -- (4.02,1.77);
\draw [dotted,color=wwttqq] (6.12,1.88) circle (0.29cm);
\draw [dotted,color=wwttqq] (2.91,2.45) circle (0.64cm);
\draw [dotted,color=wwttqq] (2.91,2.45) circle (0.41cm);
\draw [dotted,color=wwttqq] (3.55,2.43)-- (4.66,2.23);
\draw [dotted,color=wwttqq] (3.51,2.23)-- (4.42,2.07);
\draw [dotted,color=wwttqq] (9.27,2.72) circle (0.41cm);
\draw [dotted,color=wwttqq] (9.27,2.72) circle (0.62cm);
\draw [dotted,color=wwttqq] (8.65,2.66)-- (7.64,2.41);
\draw [dotted,color=wwttqq] (8.71,2.45)-- (7.98,2.26);
\draw [shift={(6.12,1.88)},color=ffzzzz]  plot[domain=2.83:3.07,variable=\t]({1*4.04*cos(\t r)+0*4.04*sin(\t r)},{0*4.04*cos(\t r)+1*4.04*sin(\t r)});
\draw [shift={(6.12,1.88)},color=wwttqq]  plot[domain=0.11:0.36,variable=\t]({1*4.01*cos(\t r)+0*4.01*sin(\t r)},{0*4.01*cos(\t r)+1*4.01*sin(\t r)});
\draw [->,color=ffzzzz] (2.08,2.21) -- (2.06,1.99);
\draw [->,color=wwttqq] (10.11,2.24) -- (10.11,2.08);
\draw [color=wwttqq] (11.04,0)-- (5.71,0);
\draw [color=wwttqq] (5.71,0)-- (5.19,-0.62);
\draw [color=wwttqq] (5.52,-0.37) -- (5.38,-0.25);
\draw [color=wwttqq] (5.19,-0.62)-- (5.33,-1.43);
\draw [color=wwttqq] (5.35,-1.01) -- (5.17,-1.04);
\draw [color=wwttqq] (5.33,-1.43)-- (6.04,-1.84);
\draw [color=wwttqq] (5.73,-1.56) -- (5.64,-1.71);
\draw [color=wwttqq] (6.04,-1.84)-- (11.04,0);
\draw [color=wwttqq] (2.86,0)-- (5.71,0);
\draw [color=wwttqq] (4.25,0.09) -- (4.25,-0.09);
\draw [color=wwttqq] (4.32,0.09) -- (4.32,-0.09);
\draw [color=wwttqq] (5.33,-1.43)-- (2.86,0);
\draw [color=wwttqq] (4.08,-0.81) -- (4.17,-0.65);
\draw [color=wwttqq] (4.02,-0.77) -- (4.11,-0.62);
\draw [dotted,color=wwttqq] (6.12,-0.88) circle (0.29cm);
\draw [dotted,color=wwttqq] (9.27,-1.72) circle (0.62cm);
\draw [dotted,color=wwttqq] (9.27,-1.72) circle (0.41cm);
\draw [dotted,color=wwttqq] (8.71,-1.45)-- (7.98,-1.26);
\draw [dotted,color=wwttqq] (8.65,-1.66)-- (7.64,-1.41);
\draw [dotted,color=wwttqq] (3.55,-1.43)-- (4.66,-1.23);
\draw [dotted,color=wwttqq] (3.51,-1.23)-- (4.42,-1.07);
\draw [dotted,color=wwttqq] (2.91,-1.45) circle (0.41cm);
\draw [dotted,color=wwttqq] (2.91,-1.45) circle (0.64cm);
\fill [color=wwttqq] (2.86,1) circle (1.0pt);
\draw[color=wwttqq] (2.72,0.71) node {$A$};
\fill [color=wwttqq] (11.04,1) circle (1.0pt);
\draw[color=wwttqq] (11.06,0.65) node {$B$};
\fill [color=wwttqq] (6.04,2.84) circle (1.0pt);
\draw[color=wwttqq] (6.2,3.05) node {$C$};
\fill [color=wwttqq] (5.71,1) circle (1.0pt);
\draw[color=wwttqq] (5.7,0.67) node {$M$};
\fill [color=wwttqq] (5.33,2.43) circle (1.0pt);
\draw[color=wwttqq] (5.34,2.71) node {$K$};
\fill [color=wwttqq] (5.19,1.62) circle (1.0pt);
\draw[color=wwttqq] (4.9,1.79) node {$L$};
\fill [color=wwttqq] (11.04,0) circle (1.0pt);
\fill [color=wwttqq] (5.71,0) circle (1.0pt);
\fill [color=wwttqq] (5.19,-0.62) circle (1.0pt);
\fill [color=wwttqq] (5.33,-1.43) circle (1.0pt);
\fill [color=wwttqq] (6.04,-1.84) circle (1.0pt);
\fill [color=wwttqq] (2.86,0) circle (1.0pt);
\end{tikzpicture}}
\caption{A scissors-shaped nice cut}
\label{scissors}
\end{figure}
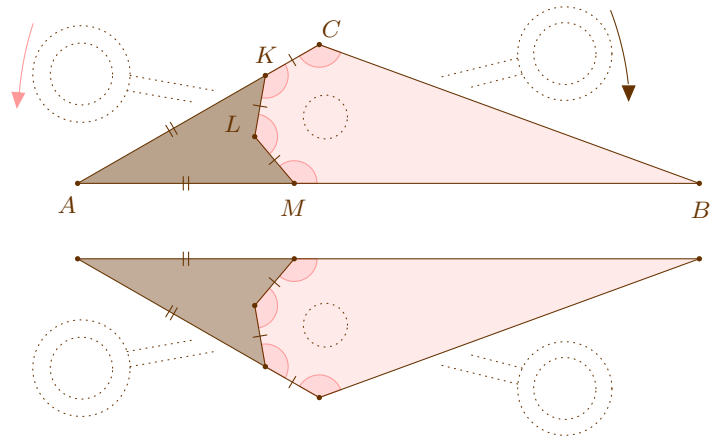

\begin{proof}[Third proof of Example~\ref{ex1}\textup{(e):} {\bf 3. A scissors-shaped nice cut}]
Take a broken line $BCKLMB$ with $CK=KL=LM$ and with equal angles $130^\circ$ between
consecutive edges, see Figure~\ref{scissors}. Then $K\in AC$, $M\in AB$ and the cut $KLM$ is nice.
Rotating the ''scissors'' $AKLM$ and $BCKLMB$ we pack the cake into the box.
\end{proof}


\section{Acknowledgements} The author is grateful to B.R.~Frenkin, A.A.~Glazyrin, I.V.~Izmestiev and M.V.~Prasolov for useful discussions. The author is also grateful to his wife Anastasia for some figures and cakes. The author was supported in part by Moebius Contest Foundation for Young Scientists and the Euler Foundation.


\begin{thebibliography}{Bol77}

\bibitem{Bol77} 
V.G. Boltianskii,
\textit{Hilbert's third problem}, Transl. by R. Silverman, V. H. Winston \& Sons, Washington D.C., 1978.




\bibitem{PSF} M. Prasolov, M. Skopenkov and B. Frenkin, Invariants of polygons, \textit{XIX Summer conference of the International mathematical tournament of towns}, 2007, \url{http://turgor.ru/lktg/2007/1/1-1en.pdf}.

\bibitem{Zak} A. Zak, Dissection of a triangle into similar triangles, \textit{Discr. Comp. Geom.} \textbf{34:2} (2005), pp. 295--312.
\end{thebibliography}
\end{document}